\documentclass[a4paper,11pt]{article}

\usepackage{lmodern}
\usepackage{enumerate}
\usepackage{amsmath,amssymb}
\usepackage{amsthm}
\usepackage{amsfonts}
\usepackage{color,soul}
\usepackage{xypic}
\usepackage{tikz}
\usepackage{accents} 
\usepackage{float} 
\usepackage{setspace} 
\usetikzlibrary{arrows, matrix, decorations.pathmorphing, backgrounds}

\newcommand{\mcal}{\mathcal}
\newcommand{\Ab}{\mathrm{Ab}}

\DeclareMathOperator{\End}{End}
\DeclareMathOperator{\Hom}{Hom}
\DeclareMathOperator{\op}{op}
\newcommand{\ModR}[1]{\operatorname{Mod}\mbox{-}{#1}}
\newcommand{\modR}[1]{\operatorname{mod}\mbox{-}{#1}}
\newcommand{\indR}[1]{\operatorname{ind}\mbox{-}{#1}}

\newcommand{\RMod}[1]{{#1}\mbox{-}\operatorname{Mod}}
\newcommand{\Rmod}[1]{{#1}\mbox{-}\operatorname{mod}}
\newcommand{\smodR}[1]{\underline{\operatorname{mod}}\mbox{-}{#1}}
\newcommand{\sModR}[1]{\underline{\operatorname{Mod}}\mbox{-}{#1}}

\DeclareMathOperator{\Zg}{Zg}
\DeclareMathOperator{\Sp}{Sp}

\theoremstyle{plain}
\newtheorem{theorem}{Theorem}[section]
\newtheorem{lemma}[theorem]{Lemma}
\newtheorem{proposition}[theorem]{Proposition}
\newtheorem{corollary}[theorem]{Corollary}

\theoremstyle{definition}

\newtheorem{example}[theorem]{Example}

\DeclareMathAlphabet{\mathantt}{OT1}{antt}{m}{n}
\def\field{\mathantt{k}}

\renewcommand{\hat}\widehat

\begin{document}
\title{Ziegler Spectra of Self Injective Algebras of Polynomial Growth}
\author{Michael Bushell}
\maketitle
\begin{abstract} 
A description of the Ziegler spectrum is given for trivial extensions of tubular algebras and related self-injective algebras of tubular type.
\end{abstract}

{
\tableofcontents
}

\section{Introduction}
An algebra is \textbf{self-injective} if it is injective as a module over itself. A great deal of work has been done classifying such algebras. A survey of some known results concerning tame self-injective algebras (and their Morita, derived, and stable, equivalence classes) can be found in \cite{skowronski}.

In particular, the self-injective algebras of polynomial growth are well understood. Many of these algebras admit universal coverings by simply-connected algebras; these coverings can be visualised as (topological) coverings of the corresponding quivers. The induced push-down functors, on finite-dimensional modules, have been used successfully in representation theory \cite{bongartz_gabriel} \cite{gabriel} \cite{dowbor_lenzing_skowronski}. Here, we extend the use of these covering techniques to determine the Ziegler spectrum of a self-injective algebra of polynomial growth.

The Ziegler spectrum of a ring $R$ is the set of isomorphism classes of indecomposable pure-injective $R$-modules, topologised as in \cite{ziegler}. For a finite-dimensional algebra, this includes all indecomposable finite-dimensional modules and, whenever the algebra is representation-infinite, there are additional infinite-dimensional points. The finite-dimensional modules are partitioned into components of the Auslander-Reiten quiver, and it is often fruitful to consider the closure of these components in the Ziegler spectrum. This is done for stable tubes in \cite{ringel3} and for ray/coray tubes in \cite{gregory}. Here we solve this problem for some non-stable tubes appearing in the AR quivers of certain self-injective algebras.

Section~\ref{sec:Preliminaries} gives a brief overview of the background material used in this paper. We go into more depth on two constructions used heavily in subsequent parts: coverings and their push-down functors; and one-point extensions/coextensions. In Section~\ref{sec:NonStableTubes} we define particular algebras containing instances of non-stable tubes and compute their Ziegler closure. In Section~\ref{sec:TrivExtCanonical} we determine the Ziegler spectrum of a trivial extension of a canonical tubular algebra. In Section~\ref{sec:ZgSelfInj} we extend this result to determine the Ziegler spectrum of the trivial extension of any tubular algebra. Further, we relate these results to the Ziegler spectrum of an arbitrary self-injective algebra of polynomial-growth. The appendix of Section~\ref{sec:Appendix} contains a number of results appearing in the literature that may be less familiar than the background material of Section~\ref{sec:Preliminaries}.

\paragraph{Acknowledgements} I would like to thank the University of Manchester and EPSRC, for their financial support, and my PhD supervisor Mike Prest, for reading drafts of this paper and giving valuable feedback.

\section{Preliminaries}\label{sec:Preliminaries}
Throughout $\field$ will denote a fixed algebraically closed field. By \textbf{algebra} we mean an associative $\field$-algebra, usually finite-dimensional, and assumed to be basic and connected. Given an algebra $R$, let $\ModR{R}$ denote the category of all right (unital) $R$-modules, let $\modR{R}$ be its full subcategory of finitely-presented modules. For any class $\mcal{C}$ of modules in $\modR{R}$, let $\mathrm{add}(\mcal{C})$ denote the closure of $\mcal{C}$ under direct summands and finite direct sums, and let $\mathrm{ind}(\mcal{C})$ denote a set of chosen representatives of each isomorphism class in $\mcal{C}$. In particular, define $\indR{R} := \mathrm{ind}(\modR{R})$.

For a finite-dimensional algebra $R$, the $\field$-space dual functor $\Hom_\field(-, \field)$, induces a duality $(-)^\star: (\modR{R})^{\op} \to \modR{R^{\op}}$ between right and left finite-dimensional modules, called the \textbf{standard duality}.

We will use the language of quivers and relations of \cite[\S2]{simson_skow1}. Every finite-dimensional (basic and connected) algebra $R$ has a presentation $R \simeq \field\mcal{Q}/I$ as a \textbf{bound quiver algebra} (i.e.~the path algebra of a quiver $\mcal{Q}$ modulo an ideal $I$ generated by a set of admissible relations). There is an equivalence $\ModR{R} \simeq \mathrm{Rep}_\field(\mcal{Q}, I)$ between the categories of $R$-modules and representations of $\mcal{Q}$ bound by $I$.

The fundamental group $\Pi(\mcal{Q}, I)$ of a quiver $\mcal{Q}$ with admissible ideal $I$ may be defined \cite{green} \cite{villa-delapena}. An algebra $R$ is \textbf{simply connected} if it is triangular (i.e. its quiver $\mcal{Q}_R$ contains no oriented cycles) and for any presentation $R \simeq \field\mcal{Q}/I$ the fundamental group $\Pi(\mcal{Q}, I)$ is trivial.

We will also use the following categorical language. If $A$ is a small $\field$-linear category, let $\RMod{A} := (A, \Ab)$ denote the category of all additive functors from $A$ to $\Ab$ (equivalent to the category of all $\field$-linear functors from $A$ to $\ModR{\field}$), similarly let $\ModR{A} := (A^{\op}, \Ab)$. We sometimes treat a functor $A \to \Ab$ using module notation as in \cite[\S 3, p.17]{mitchell}.

If $R$ is an algebra with enough local idempotents $\{ e_i \}_{i \in I}$, define a corresponding category -- also denoted $R$ -- with object set $I$ and morphism spaces $R(i, j) := e_jRe_i$ for each $i, j \in I$. Composition is given by multiplication in $R$, i.e. $g \circ f := gf$ for $f \in R(i, j)$ and $g \in R(j, k)$. Then $R$ is a $\field$-linear category (it is \textbf{locally-bounded} in the sense of \cite[\S2]{bongartz_gabriel}) and $(R^{\op}, \Ab)$ is equivalent to the category of right $R$-modules. Under this equivalence, the indecomposable projective module $e_iR$ corresponds to the representable functor $R(-, i)$ for all $i \in I$. If $R = \field\mcal{Q}/I$ and $\{e_i\}_{i \in \mcal{Q}}$ are the idempotents corresponding to the vertices of $\mcal{Q}$, then the category $R^{\op}$ is the \textbf{bound path category} of $\mcal{Q}$ modulo $I$. We will freely move between modules, functors, and quiver representations, whenever convenient.

If $F: A\to B$ is a $\field$-linear functor between $\field$-linear categories, then \textbf{restriction along $F$} is the functor $\mathrm{res}_F: (B, \Ab) \to (A, \Ab)$ defined by $M \mapsto M \circ F$. Similarly $\mathrm{res}_F: (B^{\op}, \Ab) \to (A^{\op}, \Ab)$ is defined (technically by restriction along $F^{\op}$ but usually we forgo writing ``$\op$'' here).

An \textbf{$A$-$B$-bimodule} is a module $M$ that is a \emph{left} $A$-module and \emph{right} $B$-module, such that $(am)b = a(mb)$ for all $a \in A$, $b \in B$ and $m \in M$. An additive functor $M: A \otimes_\field B^{\op} \to \Ab$ (or equivalently $M: B^{\op} \otimes_\field A \to \Ab$) defines an $A$-$B$-bimodule for which $\field$ acts centrally.

\paragraph{}
For a locally-bounded algebra $R$, the category $\modR{R}$ has Auslander-Reiten (= almost-split) sequences and a corresponding translation quiver $\Gamma_R := \Gamma(\modR{R})$. We denote by $\tau_R$ (or just $\tau$) the AR translate of $\Gamma_R$. Each connected component $\Gamma^\prime$ of $\Gamma_R$ defines a \textbf{component} $\mcal{C}_{\Gamma^\prime} := \mathrm{add}(\Gamma^\prime)$ of $\modR{R}$. If $\mcal{C}$ is such a component, then $\mcal{C}$ itself has almost-split sequences and $\Gamma(C) = \Gamma^\prime$.

A component $\mcal{C}$ is called \textbf{stable} if the translate $\tau_A$ is defined for every vertex of $\Gamma(\mcal{C})$. A component $\mcal{C}$ is called \textbf{standard} if $\mathrm{ind}(\mcal{C})$ (considered as a full subcategory of $\modR{R}$) is equivalent to the bound path category of $\Gamma(\mcal{C})$ modulo the ``mesh relations'' \cite[\S 2.3]{ringel}. If $\mcal{C}_1$ and $\mcal{C}_2$ are components, then $\mcal{C}_1 \vee \mcal{C}_2$ denotes their additive closure $\mathrm{add}(\mcal{C}_1 \cup \mcal{C}_2)$ in $\modR{R}$. Additionally, \textbf{preprojective} and \textbf{preinjective} components are defined in \cite[2.1]{ringel}; and \textbf{separating} components and \textbf{tubular families} are defined in \cite[3]{ringel}. For the definition of a \textbf{tube} see \cite[3.1]{ringel} or \cite{deste_ringel}. Roughly speaking, a tube is a component $\mcal{C}$ such that $\Gamma(\mcal{C})$ contains a cyclic path and has as a geometric realisation the infinite cylinder $S^1 \times \mathbb{R}^+_0$. The \textbf{mouth} of a tube is then defined as the vertices corresponding to the subset $S^1 \times \{ 0 \}$.

\paragraph{}
For the definition of \textbf{tame}, \textbf{wild}, \textbf{domestic}, and \textbf{polynomial growth}, for finite-dimensional algebras see \cite{crawley_boevey2} and \cite[XIX.3.6]{simson_skow3} and, for locally-bounded categories, see \cite{dowbor_skowronski}.

\paragraph{}
A \textbf{Euclidean quiver} is a quiver whose underlying graph is one of the extended Dynkin (= Euclidean) diagrams $\widetilde{\mathbb{A}}_n$ (for $n \ge 1$), $\widetilde{\mathbb{D}}_n$ (for $n \ge 4$), or $\widetilde{\mathbb{E}}_n$ (for $n=6,7,8$) \cite[VII.2]{simson_skow1}, and in the first case the quiver must not be an oriented cycle. The path algebras of Euclidean quivers are the domestic hereditary algebras \cite[3.6]{ringel}. For such an algebra $R$ we have
\begin{equation*}
\modR{R} = \mcal{P} \vee \mcal{T} \vee \mcal{Q}
\end{equation*}
where $\mcal{P}$ is a preprojective component, $\mcal{T} = \bigvee_{\lambda \in \mathbb{P}^1(\field)} \mathbb{T}_\lambda$ is a separating family of standard stable tubes called the \textbf{regular} components; and $\mcal{Q}$ is a preinjective component.

A \textbf{tilted algebra} (of Euclidean type) is an algebra of the form $\End_R(T)$ where $R$ is a domestic hereditary algebra and $T$ is a tilting $R$-module (see \cite[4.1]{ringel} for a definition of \textbf{tilting module}). Furthermore, if $T$ is preprojective (or preinjective) as an $R$-module, then the tilted algebra $\End_R(T)$ is called a \textbf{tame concealed algebra} (of Euclidean type).

A \textbf{tubular algebra} is a tubular extension of a tame concealed algebra of extension type $(2,2,2,2)$, $(3,3,3)$, $(2,4,4)$, or $(2,3,6)$ \cite[\S 5]{ringel} (see \cite[4.7]{ringel} for the definition of a \textbf{tubular extension}). For such an algebra $R$ we have
\begin{equation*}
\modR{R} = \mcal{P} \vee \bigvee_{q \in \mathbb{Q}^\infty_0} \mcal{T}_q \vee \mcal{Q}
\end{equation*}
where $\mcal{P}$ is a preprojective component; for all $q \in \mathbb{Q}^\infty_0 = \mathbb{Q} \cup \{ 0, \infty \}$, $\mcal{T}_q$ is a separating $\mathbb{P}^1(\field)$-family of standard tubes (consisting of stable tubes when $q \ne 0, \infty$); and $\mcal{Q}$ is a preinjective component \cite[5.2]{ringel}.

\paragraph{}
A exact sequence in $\ModR{R}$ is \textbf{pure-exact} if it remains exact after tensoring with any (finitely-presented) left $R$-module. For many equivalent definitions see \cite[5.2]{prest4}. A module $M \in \ModR{R}$ is \textbf{pure-injective} if any pure-exact sequence beginning at $M$ is split exact (i.e. $M$ is injective over pure embeddings).

The isomorphism classes of indecomposable pure-injective $R$-modules form a set, denoted $\Zg_R$, called the (right) \textbf{Ziegler spectrum} of $R$. This set was first topologised in model theoretic terms \cite[4.9]{ziegler} but we use the following equivalent description. A full subcategory of $\mcal{D} \subseteq \ModR{R}$ is \textbf{definable} if closed under products, direct limits, and pure submodules (i.e. images of pure embeddings). If $\mcal{D}$ is a definable subcategory of $\ModR{R}$, then $\Zg(\mcal{D}) := \mcal{D} \cap \Zg_R$ is a typical closed subset of $\Zg_R$.

A functor is an \textbf{interpretation functor} if it commutes with products and direct limits. An interpretation functor $F: \ModR{R} \to \ModR{S}$ induces a morphism from the lattice of closed sets of $\Zg_R$ to the lattice of closed sets of $\Zg_S$ \cite[15.2]{prest4}. 

Definable subcategories can be defined in terms of kernels of interpretation functors \cite[10.2.43]{prest}. Examples include the ``hom-orthogonal'' classes of finitely-presented modules, i.e. given $M \in \modR{R}$, then $\{ X \mid (M, X) = 0 \}$ and, if $R$ is an (artin) algebra, $\{ X \mid (X, M) = 0 \} = \{ X \mid X \otimes M^\star = 0 \}$ are definable subcategories of $\ModR{R}$ \cite[10.2.35, 10.2.36]{prest}.

\begin{proposition}{\cite[15.2,15.3]{prest4}}\label{prop:InterpFunc}
Let $\mcal{D} \subseteq \ModR{S}$ be a definable subcategory and $F: \mcal{D} \to \ModR{R}$ is an interpretation functor. If $X \subseteq \Zg(\mcal{D})$ is a closed subset and $F$ preserves indecomposability on $X$ (i.e. $F(M)$ is indecomposable for all $M \in X$), then $F$ induces a closed and continuous map $X \to \Zg_S$.
\end{proposition}

For example, if $f: A \to B$ is an algebra morphism, then restriction along $f$, i.e. $\mathrm{res}_f: \ModR{B} \to \ModR{A}$, is an interpretation functor. It is full and faithful if and only if $f$ is a ring epimorphism \cite[XI.1.2]{stenstrom} and, if this is the case, then we get a closed embedding $\Zg_B \to \Zg_A$.

\begin{proposition}\label{prop:ZgFdAlg}
Let $R$ be a finite-dimensional algebra. 
\begin{enumerate}[(i)]
\item Every indecomposable finite-dimensional module gives a point of $\Zg_R$ which is both open (i.e. isolated) and closed.
\item Any set of finite-dimensional points carries the discrete topology.
\item Any infinite set of finite-dimensional points contains, in its closure, at least one infinite-dimensional point and no extra finite-dimensional points.
\item The closure of the set of all finite-dimensional points equals $\Zg_R$.
\end{enumerate}
\begin{proof}
For (i), if $M$ is a finite-dimensional $R$-module, then $M$ has finite endo-length and is therefore pure-injective by \cite[4.4.24]{prest}. If $M$ is also indecomposable, then $\{ M \}$ is a closed subset of $\Zg_R$ by \cite[4.4.30]{prest} and an open subset by \cite[5.3.33]{prest}. Now (ii) is immediate from (i). For (iii), $\Zg_R$ is compact by \cite[5.1.23]{prest}, so if $X \subseteq \Zg_R$ is an infinite subset, its closure $\mathrm{cl}(X)$ contains a limit point. It follows from (ii) that, if $X$ consists of finite-dimensional points, then any limit point must be infinite-dimensional. Finally (iv) is \cite[5.3.36]{prest}.
\end{proof}
\end{proposition}

Therefore, to determine the Ziegler spectrum of a finite-dimensional algebra, it is effectively enough to provide a cover of finitely many closed subsets (whose relative topology is understood) containing almost all finite-dimensional points.

\paragraph{}
The Ziegler spectrum of a tame hereditary algebra is given in \cite{prest2} and \cite{ringel2}. Every infinite-dimensional point of $\Zg_R$, for $R$ tame hereditary and connected, lies in the closure of a tube. $\Zg_R$ contains a unique \textbf{generic} module (i.e.~an indecomposable module of infinite dimension but having finite length over its endomorphism ring \cite{crawley_boevey2}).

For information on the Ziegler spectrum of a tubular algebra see \cite{harland_prest}, \cite{reiten_ringel}, and \cite{harland}. Every (non-preprojective, non-preinjective) point of $\Zg_R$, for $R$ a tubular algebra, has a \textbf{slope} $r \in \mathbb{R}^\infty_0$, and for a finite-dimensional point this slope is equal to the index of the tubular family in which it lies. The points of slope $r$ generate a definable subcategory of $\ModR{R}$ and form a closed subset of $\Zg_R$. For rational $r \in \mathbb{Q}^\infty_0$ this closed subset corresponds to the closure of the tubes of index $r$ \cite[3]{ringel3}.

\paragraph{}
Let $Z$ be a topological space, denote by $Z^\prime \subseteq Z$ the closed subset of non-isolated points, this set is the \textbf{CB (Cantor-Bendixson) derivative} of $Z$. Set $Z^{(0)} := Z^\prime$ and $Z^{(n)} := (Z^{(n-1)})^\prime$ for all $n \ge 1$. These derivatives can be continued transfinitely \cite[\S 5.3.6]{prest} and must eventually stabilize, so there exists an ordinal $\alpha$ such that $Z^{(\beta)} = Z^{(\alpha)}$ for all $\beta > \alpha$ --- define $Z^{(\infty)} := Z^{(\alpha)}$. If there exists $\alpha$ such that $Z^{(\alpha-1)} \ne \emptyset$ and $Z^{(\alpha)} = \emptyset$, then say $Z$ has \textbf{CB rank $\alpha-1$}, denoted $\mathrm{CB}(Z) = \alpha-1$. Otherwise, say $Z$ has undefined rank, denoted $\mathrm{CB}(Z) = \infty$. For a point $p \in Z$, set $\mathrm{CB}(p) = \alpha$ if $p \in Z^{(\alpha)} \backslash Z^{(\alpha+1)}$ (i.e.~if $p$ is isolated in $Z^{(\alpha)}$) and $\mathrm{CB}(p) = \infty$ if no such $\alpha$ exists.

For the definition of the \textbf{isolation condition} on a closed subset of the Ziegler spectrum, and for what it means to be \textbf{isolated by an $M$-minimal pp-pair}, we refer to \cite[\S 5.3.1, \S 5.3.2]{prest}. The \textbf{m-dimension} $\mathrm{mdim}(\mcal{D})$ of a definable category $\mcal{D}$ is defined in \cite[\S 7.2]{prest}.

\begin{proposition}[{\cite[5.3.60]{prest}}]
If a closed subset $X \subseteq \Zg_R$ satisfies the isolation condition, then $CB(X) = \mathrm{mdim}(\mcal{X})$ where $\mcal{X}$ is the definable subcategory of $\ModR{R}$ corresponding to $X$.
\end{proposition}

\paragraph{}
We have said an algebra $A$ is \textbf{self-injective} if $A_A$ is an injective $A$-module. For algebras without unity, by self-injective we mean that the projective and injective (left or right) modules coincide. This same definition can be used for finite-dimensional algebras by \cite[IV.3.1]{auslander_smalo}.

\paragraph{}
Let $R$ be a finite-dimensional algebra and $M$ a finite-dimensional $R$-$R$-bimodule. The trivial extension $R \ltimes M$ of $R$ by $M$ is the $\field$-vector space $R \oplus M$ with pointwise addition and multiplication defined by
\begin{equation*}
(r,x)(s,y) = (rs, xs+ry)
\end{equation*}
for all $r,s \in R$ and $x, y \in M$. We may identify $M$ with the ideal $0 \oplus M$ of $R \ltimes M$, then $M^2 = 0$ and the canonical projection $R \ltimes M \to R$ has kernel $M$. Using $M = R^\star$ we call $R \ltimes R^\star$ simply the \textbf{trivial extension} of $R$ (here $R^\star$, the $\field$-dual space, has the $R$-$R$-bimodule structure inherited from $R$).

Consider $R$ as a category with object set $I$ corresponding to idempotents $\{ e_i \}_{i \in I}$ of $R$. Define the \textbf{repetitive category} $\hat{R}$ of $R$ with object set $I \times \mathbb{Z}$ and morphism spaces
\begin{equation*}
((i, m), (j, n)) = \begin{cases}
e_jRe_i & \text{if } n=m,\\
(e_i R e_j)^\star & \text{if } n=m+1,\\
0 & \text{otherwise}.
\end{cases}
\end{equation*}
Composition in $\hat{R}$ is induced by multiplication in $R$ and the $R$-$R$-bimodule structure of $R^\star$ --- note that $(e_i R e_j)^\star \simeq e_j R^\star e_i$. The category $\hat{R}$ is defined up to isomorphism \cite[1.1]{ohnuki_takeda_yamagata}. There exists an evident automorphism $\nu: \widehat{R} \to \widehat{R}$, called the \textbf{Nakayama automorphism}, defined on objects by $(i, m) \mapsto (i, m+1)$.

Both $\hat{R}$ and $R \ltimes R^\star$ are \textbf{self-injective} \cite{hughes_waschbusch}. 

\subsection{Covering and push-down functors}\label{sec:Coverings}
Let $A$ and $B$ be small $\field$-linear categories. Given a functor $F: A \to B$ and an object $b \in B$, introduce the notation $a/b$ as shorthand for $a \in F^{-1}(b)$ (i.e. $F(a) = b$).

 A $\field$-linear functor $F: A \to B$ is a \textbf{covering functor} \cite[3.1]{bongartz_gabriel} if, for all $b_1, b_2 \in B$, the following two conditions are satisfied:
\begin{enumerate}[(a)]
\item For all $a_1/b_1$ the maps $\{ F_{a_1,a_2}: A(a_1, a_2) \to B(b_1, b_2) \mid a_2/b_2 \}$ are the components of a bijection
\begin{equation}\label{eq:CovProp1}
\bigoplus_{a_2/b_2} A(a_1, a_2) \to B(b_1, b_2)
\end{equation}
\item For all $a_2/b_2$ the maps $\{ F_{a_1,a_2}: A(a_1, a_2) \to B(b_1, b_2) \mid a_1/b_1 \}$ are the components of a bijection
\begin{equation}\label{eq:CovProp2}
\bigoplus_{a_1/b_1} A(a_1, a_2) \to B(b_1, b_2)
\end{equation}
\end{enumerate}

For a covering functor $F: A \to B$ it has become usual to call restriction along $F$ (i.e. the functor $\mathrm{res}_F: \RMod{B} \to \RMod{A}$) the \textbf{pull-up} functor and its left adjoint -- denoted $F_\lambda: \RMod{A} \to \RMod{B}$ -- the \textbf{push-down} functor. The push-down functor $F_\lambda$ is defined explictly in \cite[3.2]{bongartz_gabriel} and we derive this definition below.

By Proposition~\ref{prop:KanExt}, $F_\lambda$ is defined for left (resp. right) $A$-modules by the $B$-$A$-bimodule $B(F-, -)$ (resp.~the $A$-$B$-bimodule $B(-, F-)$). The bijection \eqref{eq:CovProp1} above is easily seen to be natural in $a_1$ and gives, for each $b_2 \in B$, an $A$-linear isomorphism
\begin{equation}\label{eq:CovProp3}
\bigoplus_{a_2/b_2} A(-, a_2) \to B(F(-), b_2)
\end{equation}
Note the right-hand side is just $\mathrm{res}_F B(-, b_2)$. In this way, each representable (right) $B$-module ``lifts'' (or ``pulls-up'') to the coproduct of all representable $A$-modules ``lying above'', in the sense that
\begin{equation*}
F_\lambda A(-, a_2) = A(-, a_2) \otimes_A B(-, F-) \simeq B(-, F(a_2)) = B(-, b_2)
\end{equation*}
for all $a_2/b_2$. The bijection \eqref{eq:CovProp2} gives the analogous result for left modules. Now the right-hand side of \eqref{eq:CovProp3} is functorial in $b_2$, so there is a unique way to make the left-hand side functorial in $b_2$ such that \eqref{eq:CovProp3} becomes a $B$-$A$-bimodule isomorphism (see Lemma~\ref{lemma:BABiMod} below).

\paragraph{}
\emph{Henceforth, fix a covering functor $F: A \to B$}. Given $\beta: b_1 \to b_2$ in $B$, for each $a_1/b_1$ there exists, according to $(a)$, a unique set $\{ \beta^{a_1}_{(a_2)}: a_1 \to a_2 \mid a_2/b_2 \}$ of morphisms in $A$ (almost all zero) such that $\beta = \sum_{a_2/b_2} F(\beta^{a_1}_{(a_2)})$. The notation $\beta^{a_1}_{(a_2)}$ is to suggest the domain $a_1$ is fixed and the codomain $a_2$ varies as we ``lift'' $\beta$ to $A$ via $F$. Similarly, for each $a_2/b_2$, property (b) gives a unique set $\{ \beta^{(a_1)}_{a_2}: a_1 \to a_2 \mid a_1/b_1 \}$ of morphisms in $A$ (almost all zero) such that $\beta = \sum_{a_1/b_1} F(\beta^{(a_1)}_{a_2})$. We are careful to distinguish between $\beta^{(a_1)}_{a_2}$ and $\beta^{a_1}_{(a_2)}$.

\begin{lemma}\label{lemma:PushDownCalc}
Given $\beta: b_1 \to b_2$, $\delta: b_2 \to b_3$ in $B$, then
\begin{equation*}
(\delta\beta)^{a_1}_{(a_3)} = \sum_{a_2/b_2} \delta^{a_2}_{(a_3)}\beta^{a_1}_{(a_2)}
\end{equation*}
for all  $a_1/b_1$ and $a_3/b_3$ in $A$.
\begin{proof}
We have
\begin{align*}
\sum_{a_3/b_3} F\left( \sum_{a_2/b_2} \delta^{a_2}_{(a_3)} \beta^{a_1}_{(a_2)} \right)
&= \sum_{a_3/b_3} \sum_{a_ 2/b_2} F(\delta^{a_2}_{(a_3)}) F(\beta^{a_1}_{(a_2)}) \\
&= \sum_{a_2/b_2} \left(\sum_{a_3/b_3} F(\delta^{a_2}_{(a_3)})\right) F(\beta^{a_1}_{(a_2)}) \\
& = \sum_{a_2/b_2} \delta F(\beta^{a_1}_{(a_2)}) \\
&= \delta\beta
\end{align*}
and equality follows by the unique lifting property.
\end{proof}
\end{lemma}

For $a \in A$ and $b \in B$ define $\Theta(a, b) := \bigoplus_{x/b} A(a, x)$. We make $\Theta$ a $B$-$A$-bimodule as follows. Note $\Theta(a,b)$ is already functorial in $a$, being a coproduct of representable functors.

For $\beta: b_1 \to b_2$ in $B$ and $x \in A$ define
\begin{equation*}
\Theta(x, \beta): \bigoplus_{a_1/b_1} A(x, a_1) \to \bigoplus_{a_2/b_2} A(x, a_2)
\end{equation*}
by
\begin{equation*}
\Theta(x, \beta) ( (\alpha_{a_1})_{a_1/b_1} ) := \left( \sum_{a_1/b_1} \beta^{a_1}_{(a_2)} \alpha_{a_1} \right)_{a_2/b_2}
\end{equation*}
for  $(\alpha_{a_1})_{a_1/b_1} \in \bigoplus_{a_1/b_1} A(x, a_1)$.

\begin{lemma}
For each $x \in A$, $\Theta(x, -): B \to \Ab$ is a well-defined functor.
\begin{proof}
The claim is that $\Theta(x, b)$ is functorial in $b \in B$. Given $a/b$ the equality $F(1_a) = 1_b = \sum_{y/b} F((1_b)^{a}_{(y)})$ implies $(1_b)^a_{(a)} = 1_a$ and $(1_b)^a_{(y)} = 0$ for $y \ne a$. It follows that $\Theta(x, 1_b) = 1_{\Theta(x, b)}$. Now given $\beta: b_1 \to b_2$ and $\delta: b_2 \to b_3$ in $B$, we have
\begin{equation*}
\Theta(x, \delta\beta)((\alpha_{a_1})_{a_1/b_1}) = \left( \sum_{a_1/b_1} (\delta\beta)^{a_1}_{(a_3)} \alpha_{a_1} \right)_{a_3/b_3}
\end{equation*}
whereas
\begin{align*}
\Theta(x, \delta)(\Theta(x, \beta)((\alpha_{a_1})_{a_1/b_1}))
&= \left( \sum_{a_2/b_2} \delta^{a_2}_{(a_3)} \left( \sum_{a_1/b_1} \beta^{a_1}_{(a_2)} \alpha_{a_1} \right) \right)_{a_3/b_3} \\
&= \left( \sum_{a_1/b_1} \left( \sum_{a_2/b_2} \delta^{a_2}_{(a_3)} \beta^{a_1}_{(a_2)} \right) \alpha_{a_1} \right)_{a_3/b_3}
\end{align*}
so $\Theta(x, \delta\beta) = \Theta(x, \delta)\Theta(x, \beta)$ by Lemma~\ref{lemma:PushDownCalc}.
\end{proof}
\end{lemma}


\begin{lemma}\label{lemma:BABiMod}
$\Theta$ is a $B$-$A$-bimodule and the bijections in \eqref{eq:CovProp1} give a bimodule isomorphism $\Theta(?, -) \simeq B(F(?), -)$.
\begin{proof}
For $b \in B$ and $\alpha: a_1 \to a_2$ in $A$, $\Theta(\alpha, b) = \bigoplus_{x/b} A(\alpha, x)$ is just the diagonal action $(\alpha_x)_{x/b} \mapsto ( \alpha_x \alpha)_{x/b}$. It is easily checked that each $\Theta(\alpha, b)$ is $B$-linear and that each $\Theta(a, \beta)$ is $A$-linear. For $a \in A$ and $b \in B$ the bijection
\begin{equation*}
\phi(a, b): \bigoplus_{x/b} A(a, x) \to B(F(a), b)
\end{equation*}
given in \eqref{eq:CovProp1} is defined by
\begin{equation*}
(\alpha_x)_{x/b} \mapsto \sum_{x/b} F(\alpha_x)
\end{equation*}
for $(\alpha_x)_{x/b} \bigoplus_{x/b} A(a, x)$.
We claimed in \eqref{eq:CovProp3} this bijection is natural in $a$, indeed, given $\alpha: a_1 \to a_2$ in $A$, $b \in B$, and $(\alpha_x)_{x/b} \in \bigoplus_{x/b} A(a_2, x)$, the equality
\begin{equation*}
\left( \sum_{x/b} F(\alpha_x) \right) F(\alpha) = \sum_{x/b} F(\alpha_x \alpha)
\end{equation*}
gives $B(F(\alpha), b) \circ \phi(a_2, b) = \phi(a_1, b) \circ \Theta(\alpha, b)$. Similarly, to prove naturality in $b$, take $a \in A$, $\beta: b_1 \to b_2$ in $B$, and $(\alpha_{a_1})_{a_1/b_1} \in \bigoplus_{a_1/b_1} A(x, a_1)$, then the equality
\begin{equation*}
\sum_{a_2/b_2} F\left( \sum_{a_1/b_1} \beta^{a_1}_{(a_2)} \alpha_{a_1} \right) = \sum_{a_1/b_1} \left( \sum_{a_2/b_2} F(\beta^{a_1}_{(a_2)}) \right) F(\alpha_{a_1}) = \beta \sum_{a_1/b_1} F(\alpha_{a_1})
\end{equation*}
gives $B(F(x), \beta) \circ \phi(a, b_1) = \phi(a, b_2) \circ \Theta(x, \beta)$, as required.
\end{proof}
\end{lemma}

\begin{proposition}
Let $F: A \to B$ be a covering functor and $F_\lambda: \RMod{A} \to \RMod{B}$ the corresponding push-down functor. Given $M \in \RMod{A}$, then $F_\lambda M$ may be defined by
\begin{equation*}
(F_\lambda M)(b) := \bigoplus_{a/b} M(a)
\end{equation*}
for $b \in B$. Such that, for $\beta: b_1 \to b_2$, $(F_\lambda M)(\beta)$ is given by
\begin{equation}\label{eq:PushDownAction}
(m_{a_1})_{a_1/b_1} \mapsto \left( \sum_{a_1/b_1} M(\beta^{a_1}_{(a_2)})(m_{a_1}) \right)_{a_2/b_2}
\end{equation}
for all $(m_{a_1})_{a_1/b_1} \in \bigoplus_{a_1/b_1} M(a_1)$.
\begin{proof}
By Lemma~\ref{lemma:BABiMod} and Proposition~\ref{prop:KanExt}, we can define $F_\lambda$ using the bimodule $\Theta(?, -) = \bigoplus_{a/(-)} A(?, a)$. Then, pointwise, for $b \in B$ we have isomorphisms
\begin{equation*}
\bigoplus_{a/b} A(-, a) \otimes_A M \xrightarrow{\sim} \bigoplus_{a/b} M(a)
\end{equation*}
Hence, for $\beta: b_1 \to b_1$, we find \eqref{eq:PushDownAction} by chasing around the following diagram
\begin{equation*}
\xymatrix{
\bigoplus_{a_1/b_1} A(-, a_1) \otimes_A M \ar[d]^{\Theta(-, \beta) \otimes_A M} & \bigoplus_{a_1/b_1} M(a_1) \ar@{-->}[d] \ar[l] \\
\bigoplus_{a_2/b_2} A(-, a_2) \otimes_A M \ar[r] & \bigoplus_{a_2/b_2} M(a_2)
}
\end{equation*}
Thus the definition of $F_\lambda$ may be taken as claimed.
\end{proof}
\end{proposition}

\begin{corollary}\label{cor:FinitePushDown}
The push-down functor $F_\lambda: \RMod{A} \to \RMod{B}$ is an interpretation if and only if $F^{-1}(b)$ is finite for all $b \in B$.
\begin{proof}
 By Proposition~\ref{prop:KanExt} the push-down functor $F_\lambda$ is an interpretation functor if and only 
\begin{equation*}
\mathrm{res}_F B(-, b) = B(F-, b) = \bigoplus_{a \in F^{-1}(b)} A(-, a)
\end{equation*}
is finitely-presented -- equivalently $F^{-1}(b)$ is finite -- for all $b \in B$.
\end{proof}
\end{corollary}

\subsubsection{Galois coverings}
Let $A$ be a locally bounded category, say a group $G$ of $\field$-linear automorphisms of $A$ is \textbf{admissible} provided $G$ acts freely on $A$ (i.e.~given $g \in G$, if $g \ne 1$, then $g \cdot a\ne a$ for all $a \in A$) and has finitely many orbits.

Given $A$ and such an admissible group $G$, the \textbf{orbit category} $A/G$ is defined with objects being the $G$-orbits of objects of $A$. A map $f: x \to y$ in $A/G$ is a family $({}_b f_a) \in \prod_{(a,b) \in x \times y} A(a, b)$ such that $g \cdot ({}_b f_a) = {}_{g \cdot b}f_{g \cdot a}$ for all $g \in G$. See \cite[\S 3]{gabriel} for further details. There is a projection functor $F: A \to A/G$ such that $F(a) = G \cdot a$.

\begin{proposition}[{\cite[3.1]{gabriel}}]\label{prop:GaloisCov}
If $A$ is a locally bounded category and $G$ an admissible group of automorphisms of $A$, then the projection $F: A \to A/G$ is a covering functor. Furthermore, $F$ is universal with respect to being $G$-invariant, i.e. $Fg = F$ for all $g \in G$ and if $E: A \to B$ is any $\field$-linear functor with this property, then $E = E^\prime F$ for a unique $E^\prime: A/G \to B$.
\end{proposition}

Covering functors of this form are called \textbf{Galois coverings}.

\begin{example}\label{ex:GaloisCoverings}
For a finite-dimensional algebra $R$, the canonical projection $\hat{R} \to R \ltimes R^\star$ is a Galois covering functor, for there exists an isomorphism $R \ltimes R^\star \simeq \hat{R}/\langle\nu\rangle$, where $\nu$ is the Nakayama automorphism of $\hat{R}$ and $\langle\nu\rangle$ is the infinite cyclic group generated by $\nu$ \cite[2.2]{hughes_waschbusch}.
\end{example}

\emph{Henceforth, fix a Galois covering $F: A \to B$} (with $B = A/G$), let $F_\lambda: \RMod{A} \to \RMod{B}$  be the corresponding push-down functor. For $X \in \RMod{B}$ recall $\mathrm{res}_F (X) = XF$. For $M \in \RMod{A}$ let us use the abbreviation $M_\lambda := F_\lambda(M)$ where convenient. As $F_\lambda$ is left adjoint to $\mathrm{res}_F$ we have an adjunction isomorphism $(M_\lambda, X) \simeq (M, XF)$ which we now explore.


For $M \in \RMod{A}$ and $a \in A$ observe that
\begin{equation*}
(M_\lambda F)(a) = \bigoplus_{a^\prime / F(a)} M(a^\prime)
\end{equation*}
For $X \in \RMod{B}$ and $b \in B$ observe that
\begin{equation*}
(XF)_\lambda(b) = \bigoplus_{a/b} XF(a) = \bigoplus_{a/b} X(b)
\end{equation*}
this is a direct sum of $|F^{-1}(b)| = |G|$ copies of $X(b)$.

Define $(\mu_M)_a: M(a) \to (M_\lambda F)(a)$ to be the canonical inclusion map and define $(\epsilon_X)_b: (XF)_\lambda(b) \to X(b)$ to be the summation map $(x_a)_{a/b} \mapsto \sum_{a/b} x_a$.

\begin{lemma}\label{lemma:GalCovAdj}
$\mu_M$ (resp.~$\epsilon_X$) defines the unit (resp.~counit) of the adjunction $F_\lambda \dashv \mathrm{res}_F$.
\begin{proof}
We must show that $\mu_M: M \to M_\lambda F$ is a well-defined natural transformation. Let $\alpha: a_1 \to a_2$ in $A$ be given and take $m \in M(a_1)$, write $(\mu_M)_{a_1}(m) = (m_{a^\prime_1})_{a^\prime_1/F(a_1)}$ (so that $m_{a^\prime_1} = 0$ if $a^\prime_1 \ne a_1$, and $m_{a_1} = m$), then
\begin{align*}
(M_\lambda F)(\alpha)( (\mu_M)_{a_1}(m) )
&= \left( \sum_{a^\prime_1/ F(a_1)} M(F(\alpha)^{a^\prime_1}_{(a^\prime_2)})( m_{a^\prime_1})  \right)_{a^\prime_2/F(a_2)} \\
&= (\mu_{M})_{a_2}( M(\alpha)(m) ))
\end{align*}
since $F(\alpha)^{a_1}_{(a^\prime_2)} = 0$ if $a^\prime_2 \ne a_2$, and $F(\alpha)^{a_1}_{(a_2)} = \alpha$, by the unique lifting property. This proves $\mu_M$ is well-defined.

To show $\epsilon_X: (XF)_\lambda \to X$ is well-defined, let $\beta: b_1 \to b_2$ be given and take $(x_{a_1})_{a_1/b_1} \in (XF)_\lambda (b_1)$, then
\begin{align*}
(\epsilon_X)_{b_2} ( (XF)_\lambda(\beta) ( (x_{a_1})_{a_1/b_1})) 
&= \sum_{a_2/b_2} \left( \sum_{a_1/b_1} (XF)(\beta^{a_1}_{(a_2)}) (x_{a_1}) \right) \\
&= \sum_{a_1/b_1} X\left( \sum_{a_2/b_2} F(\beta^{a_1}_{(a_2)}) \right) (x_{a_1}) \\
&= \sum_{a_1/b_1} X(\beta)(x_{a_1}) \\
&= X(\beta)(\sum_{a_1/b_1} x_{a_1}) \\
&= X(\beta)( (\epsilon_X)_{b_1} ((x_{a_1})_{a_1/b_2}) )
\end{align*}
as required.

It is easily verified that $\mu_M$ is natural in $M$ and $\epsilon_X$ is natural in $X$. By \cite[IV.1.Th.~2]{maclane} the proof is complete if we establish the following ``triangle'' identities:
\begin{enumerate}[(a)]
\item $(\epsilon_X)_{F(a)} \circ (\mu_{XF})_a = 1_{XF(a)}$ for all $X \in \RMod{B}$ and $a \in A$,
\item $(\epsilon_{M_\lambda})_b \circ (F_\lambda(\mu_M))_b = 1_{M_\lambda(b)}$ for all $M \in \RMod{A}$ and $b \in B$.
\end{enumerate}
These are easily verified, for each instance of (a) or (b) amounts to an inclusion followed by a summation, which, by definition, is the necessarily identity morphism.
\end{proof}
\end{lemma}

For $g \in G$ and $M \in \RMod{A}$ define ${}^g M = Mg^{-1}$. One can show $M_\lambda \simeq ({}^g M)_\lambda$ for all $g \in G$ and that $M_\lambda F \simeq \bigoplus_{g \in G} {}^g M$ \cite[3.2]{gabriel}. Under this isomorphism, the unit $\mu_M$ can be expressed as the canonical inclusion map $\mu_M: M \to \bigoplus_{g \in G} {}^g M$.

\begin{lemma}\label{lemma:PushDownBij}
If $M, N \in \RMod{A}$ are such that $(M, {}^g N) = 0$ for all $g \ne 1$, then $F_\lambda: (M, N) \to (M_\lambda, N_\lambda)$ is a bijection.
\begin{proof}
We have an injection $t: \bigoplus_{g \in G} (M, {}^g N) \to (M, \bigoplus_{g \in G} {}^g N)$ given by $t( (f_g)_{g \in G} ) \mapsto \sum_{g \in G} u_g f_g$ where each $u_h: {}^h N \to \bigoplus_{g \in G} {}^g N$ is a canonical inclusion \cite[Prop.~5.1]{popescu}. The image of $t$ consists of all $f: M \to \bigoplus_{g \in G} {}^g N$ such that $p_h f = 0$ for almost all $h \in G$, where $p_h: \bigoplus_{g \in G} {}^g N \to {}^h N$ is the canonical projection. In particular, if $(M, {}^g N) = 0$ for almost all $g \in G$, then $t$ is a bijection; this certainly holds by our assumption. Note also that $u_1: N \to \bigoplus_{g \in G} {}^g N$ is the unit $\mu_N: N \to \bigoplus_{g \in G} {}^g N$.

This implies that $f \mapsto \mu_N f$ defines a bijection $(M, N) \to (M, N_\lambda F)$. Now composing with the adjunction isomorphism $(M, N_\lambda F) \to (M_\lambda, N_\lambda)$ gives a bijection $(M, N) \to (M_\lambda, N_\lambda)$ defined by
\begin{equation*}
f \mapsto \epsilon_{N_\lambda} \, F_\lambda(\mu_N f) = \epsilon_{N_\lambda} F_\lambda(\mu_N) F_\lambda(f) = F_\lambda(f)
\end{equation*}
using the triangle identity $\epsilon_{N_\lambda} F_\lambda(\mu_N) = 1_{N_\lambda}$.
\end{proof}
\end{lemma}

\begin{lemma}\label{lemma:GalProp1}
For all $\beta: b_1 \to b_2$ in $B$, $a_1/b_1$, $a_2/b_2$, we have $\beta^{a_1}_{(a_2)} = \beta^{(a_1)}_{a_2}$.
\begin{proof}
Let $\beta: b_1 \to b_2$ in $B$ and $g \in G$ be given and fix $a_1/b_1$, then
\begin{equation*}
\sum_{a_2/b_2} F(\beta^{a_1}_{(a_2)}) = \beta = \sum_{a_2/b_2} F(g \cdot \beta^{a_1}_{(a_2)})
\end{equation*}
by $G$-invariance of $F$. Now $g \cdot \beta^{a_1}_{(a_2)}$ is a morphism $g \cdot a_1 \to g \cdot a_2$ in $A$, hence $g \cdot \beta^{a_1}_{(a_2)} = \beta^{g \cdot a_1}_{(g \cdot a_2)}$ by uniqueness. Similarly $g \cdot \beta^{(a_1)}_{a_2} = \beta^{(g \cdot a_1)}_{g \cdot a_2}$. Now, fixing $a_1/b_1$ and $a_2/b_2$ we have

\begin{equation*}
\sum_{g \in G} F(\beta^{(g \cdot a_1)}_{a_2}) = \beta = \sum_{g \in G} F(\beta^{a_1}_{(g^{-1} \cdot a_2)}) = \sum_{g \in G} F(g \cdot \beta^{a_1}_{(g^{-1} \cdot a_2)}) = \sum_{g \in G} F(\beta^{g \cdot a_1}_{(a_2)})
\end{equation*}
from which it follows that $\beta^{a_1}_{(a_2)} = \beta^{(a_1)}_{a_2}$ by uniqueness.
\end{proof}
\end{lemma}

Let $F_\rho: \RMod{A} \to \RMod{B}$ denote the right adjoint to $\mathrm{res}_F$. $F_\rho$ can be derived similarly to $F_\lambda$. For $M \in \RMod{R}$ we find that
\begin{equation*}
(F_\rho M)(b) := \prod_{a/b} M(a)
\end{equation*}
for $b \in B$, and for $\beta: b_1 \to b_2$, then $(F_\rho M)(\beta)$ is defined by
\begin{equation}\label{eq:PushDownAction2}
(m_{a_1})_{a_1/b_1} \mapsto \left( \sum_{a_1/b_1} M(\beta^{(a_1)}_{a_2})(m_{a_1}) \right)_{a_2/b_2}
\end{equation}
for all $(m_{a_1})_{a_1/b_1} \in \prod_{a_1/b_1} M(a_1)$.

\begin{corollary}\label{cor:FlSubFunc}
If $F$ is a Galois covering, then $F_\lambda$ is a subfunctor of $F_\rho$.
\begin{proof}
The canonical monomorphism
\begin{equation}\label{eq:FlSubFunc}
F_\lambda M (b) = \bigoplus_{a/b} M(a) \to \prod_{a/b} M(a) = F_\rho M(b)
\end{equation}
is natural in $b \in B$ by Lemma~\ref{lemma:GalProp1} (since then \eqref{eq:PushDownAction} and \eqref{eq:PushDownAction2} coincide).
\end{proof}
\end{corollary}

Given a collection of $A$-modules $\mcal{M} = \{ M_i \}_{i \in I}$ define
\begin{equation*}
\mathrm{Supp}_A(\mcal{M}) := \bigcup_{i \in I} \{ a \in A \mid M_i(a) \ne 0 \}
\end{equation*}
and say $\mcal{M}$ is \textbf{finitely-supported} if $\mathrm{Supp}_A(\mcal{M})$ is finite.

\begin{lemma}\label{lemma:CovFinSupp}
If $M \in \RMod{A}$ is finitely-supported, then $F_\lambda M = F_\rho M$.
\begin{proof}
This is immediate from (the proof of) Corollary~\ref{cor:FlSubFunc}, for \eqref{eq:FlSubFunc} is always an isomorphism (by assumption, the product is finite and coincides with the direct sum).
\end{proof}
\end{lemma}

\begin{corollary}\label{cor:PushDownLimits}
Let $F$ be a Galois covering and $F_\lambda: \RMod{A} \to \RMod{B}$ the corresponding push-down functor. If $\mcal{M} = \left(\{ M_i \}_{i \in I}, (\gamma_{i,j}: M_j \to M_i)_{i \le j} \right)$ is any (inverse) system in $\RMod{A}$ such that $\{ M_i \}_{i \in I}$ is a finitely-supported set of $A$-modules, then $F_\lambda \varprojlim M_i = \varprojlim F_\lambda M_i$.
\begin{proof}
The set $\{ M_i \}_{i \in I}$ being finitely-supported implies $\varprojlim_{i \in I} M_i$ is finitely-supported and each $M_i$ is finitely-supported, hence
\begin{equation*}
F_\lambda\varprojlim M_i = F_\rho\varprojlim M_i = \varprojlim F_\rho M_i = \varprojlim F_\lambda M_i
\end{equation*}
by Lemma~\ref{lemma:CovFinSupp} and the fact that $F_\rho$ (being a right adjoint) commutes with limits.
\end{proof}
\end{corollary}

%

\begin{corollary}\label{cor:DefInterp}
Let $F: A \to B$ be a Galois covering. If $\mcal{D}$ a definable subcategory of $\RMod{A}$ whose class of modules is finitely-supported, then the restriction $F_\lambda|_{\mcal{D}}: \mcal{D} \to \RMod{B}$, of the push-down functor, is an interpretation functor.
\begin{proof}
As a left-adjoint, the push-down functor preserves direct limits. By Corollary~\ref{cor:PushDownLimits}, the restriction to $\mcal{D}$ then also commutes with products, and is therefore an interpretation functor.
\end{proof}
\end{corollary}

\paragraph{}
A subcategory $A^\prime \subseteq A$ is \textbf{convex} if $A^\prime$ is closed under factorisation of morphisms, i.e. for any morphism $\alpha: x \to z$ in $A^\prime$, if $\alpha = \alpha_2\alpha_1$ for some $\alpha_1: x \to y$ and $\alpha_2: y \to z$ in $A$, then $y$, $\alpha_1$, and $\alpha_2$ are in $A^\prime$. If $A^\prime \subseteq A$ is a full convex subcategory, define
\begin{equation*}
I(x,y) = \begin{cases}
A(x,y) & \text{if } x \not\in A^\prime \text{ or } y \not\in A^\prime,\\
0 & \text{otherwise}.
\end{cases}
\end{equation*}
Then $I$ is an ideal of $A$ (it is the ideal generated by the set of objects $A \backslash A^\prime$) and the quotient category $A/I$ is defined with the same objects as $A$ and with hom-spaces $(A/I)(x,y) := A(x,y)/I(x,y)$ for all $x, y \in A$. Let $A^\prime_0$ be the category obtained from $A^\prime$ by adding a null object $0$ (i.e. $0$ is both initial and terminal in $A^\prime_0$, thus all hom-spaces to and from $0$ are zero). Now the function
\begin{equation*}
a \mapsto \begin{cases}
a & \text{if } a \in A^\prime,\\
0 & \text{if } a \not\in A^\prime,
\end{cases}
\end{equation*}
for $a \in A$, then extends to a functor $A/I \to A^\prime_0$ which is easily seen to be an equivalence. Composing with the canonical projection $A \to A/I$ gives a functor $\pi: A \to A^\prime_0$ which is full and surjective on objects, and therefore $\mathrm{res}_\pi: \RMod{A^\prime_0} \to \RMod{A}$ is a full and faithful embedding. We have a clear equivalence $\ModR{A^\prime} \simeq \ModR{A^\prime_0}$ and, via the composition with $\mathrm{res}_\pi$, we may consider $\ModR{A^\prime}$ as a definable subcategory of $\ModR{A}$. 

\begin{corollary}\label{cor:ConvexRes}
Let $F: A \to B$ be a Galois covering. If $A^\prime \subseteq A$ is a finite full and convex subcategory of $A$, then the restriction $F_\lambda|_{\RMod{A^\prime}}: \RMod{A^\prime} \to \RMod{B}$, of the push-down functor, is an interpretation functor. Futhermore, if $A^\prime$ intersects each $G$-orbit of $A$ at most once, then the restriction $F_\lambda|_{\RMod{A^\prime}}$ is full and faithful.
\begin{proof}
Apply Corollary~\ref{cor:DefInterp} with $\mcal{D} = \RMod{A^\prime}$. For the final statement, if $A^\prime$ intersects each $G$-orbit of $A$ at most once, then for all $M, N \in \RMod{A^\prime}$ and $g \in G \backslash \{ 1 \}$, the assumption implies $\mathrm{Supp}({}^g N)$ lies outside $A^\prime$, therefore $\mathrm{Supp}(M) \cap \mathrm{Supp}({}^g N) = \emptyset$ and $(M, {}^g N) = 0$. Thus Lemma~\ref{lemma:PushDownBij} applies and shows $F_\lambda|_{\RMod{A^\prime}}$ is full and faithful. 
\end{proof}
\end{corollary}

\subsection{One-point extensions and coextensions}\label{sec:OnePointExt}
Let $R$ be a finite-dimensional algebra and $X \in \modR{R}$ a finite-dimensional \emph{right} $R$-module. The \textbf{one-point extension of $R$ by $X$} \cite[2.2.5]{ringel} is the following algebra of matrices:
\begin{equation*}
R[X] := \begin{pmatrix} R & 0 \\ X & \field \end{pmatrix} = \left\{\begin{pmatrix} r & 0 \\ x & \lambda \end{pmatrix} \mid r \in R, x \in X, \lambda \in \field \right\}
\end{equation*}
The \textbf{one-point coextension of $R$ by $X$} is $[X]R := (R^{\op}[X^\star])^{\op}$ and is isomorphic to the following algebra of matrices:
\begin{equation*}
[X]R \simeq \begin{pmatrix} \field & 0 \\ X^\star & R \end{pmatrix}
\end{equation*}
In these definitions, we are using the $\field$-$R$-bimodule (resp. $R$-$\field$-bimodule) structure of $X$ (resp. $X^\star$)  induced by the commutativity of $\field$. Note, we have the standard duality
\begin{equation*}
\modR{[X]R} \simeq (\modR{R^{\op}[X^\star]})^{\op}
\end{equation*}
and results concerning finite-dimensional modules over one-point extensions dualise to one-point coextensions. For more on one-point extensions of bound quiver algebras in particular see \cite[XV.1]{simson_skow3}.

\begin{example}\label{ex:CanonicalTubular}
Let $A$ be the path algebra of the following quiver (below left) and $X$ the module (below middle).
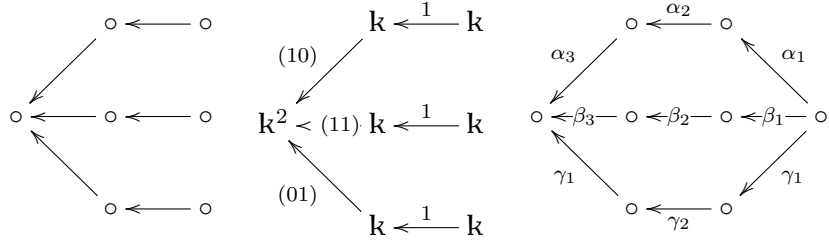
\begin{equation*}
\xymatrix{
& \circ \ar[dl] & \circ \ar[l] \\
\circ & \circ \ar[l] & \circ \ar[l] \\ 
& \circ \ar[ul] & \circ \ar[l]
} \quad
\xymatrix{
& \field \ar[dl]_{(1 0)} & \field \ar[l]_{1} \\
\field^2 & \field \ar[l]|{\quad(1 1)} & \field \ar[l]_{1} \\ 
& \field \ar[ul]^{(0 1)} & \field \ar[l]_{1}
} \quad
\xymatrix{
& \circ \ar[dl]_{\alpha_3} & \circ \ar[l]_{\alpha_2} \\
\circ & \circ \ar[l]|{\beta_3} & \circ \ar[l]|{\beta_2} & \circ \ar[l]|{\beta_1} \ar[ul]_{\alpha_1} \ar[dl]^{\gamma_1} \\ 
& \circ \ar[ul]^{\gamma_1} & \circ \ar[l]^{\gamma_2}
}
\end{equation*}
The one-point extension $A[X]$ is the \textbf{canonical algebra} of type $(3,3,3)$ defined in \cite[3.3.7]{ringel}. It is given by the quiver (above right) with the relation $\alpha_1\alpha_2\alpha_3 + \beta_1\beta_2\beta_3 + \gamma_1\gamma_2\gamma_3 = 0$.
\end{example}

If $R$ has a complete set of local idempotents $\{ e_1, \dots, e_n \}$, then $R[X]$ has a complete set of local idempotents
\begin{equation*}
\left\{ f_i = \begin{pmatrix} e_i & 0 \\ 0 & 0 \end{pmatrix} \text{ for } i=1,\dots,n, \text{ and } f_\omega = \begin{pmatrix} 0 & 0 \\ 0 & 1 \end{pmatrix} \right\}
\end{equation*}
and we call $\omega$ the \textbf{extension vertex} of $R[X]$. The ideal $I = f_\omega R[X]$ is such that $R \simeq R[X]/I$ and the canonical projection $\pi: R[X] \to R$ induces a full and faithful embedding $\mathrm{res}_\pi: \ModR{R} \to \ModR{R[X]}$ which we call the \textbf{zero embedding}, since it is given by ``extending by 0''. 

Consider $R$ (resp.~$R[X]$) as categories with objects $\{ 1, \dots, n \}$ (resp. $\{ 1, \ldots, n, \omega \}$). We regard $R$ as a full subcategory of $R[X]$ by identifying the hom-spaces
\begin{equation}\label{eq:RXHoms}
R[X](i, j) = f_jR[X]f_i \simeq e_jRe_i = R(i,j)
\end{equation}
for all $i, j \in \{ 1, \dots, n \}$. Let $F: R \to R[X]$ denote the inclusion functor and $\mathrm{res}_F: \ModR{R[X]} \to \ModR{R}$ denote restriction along $F$ (stricly speaking, along $F^{\op}$), and simlarly $\mathrm{res}_F: \RMod{R[X]} \to \RMod{R}$.

\begin{lemma}\label{lemma:RXProjRes}
For all $i \in \{ 1, \dots, n, \omega \}$ we have
\begin{align*}
\mathrm{res}_F R[X](-, i) &\simeq \begin{cases}
R(-, i) &\text{if } i \ne \omega,\\
X &\text{if } i = \omega.
\end{cases} \\
\mathrm{res}_F R[X](i,-) &\simeq \begin{cases}
R(i,-) &\text{if } i \ne \omega,\\
0 &\text{if } i = \omega.
\end{cases}
\end{align*}
as $R$-modules.
\begin{proof}
Pointwise, this is clear from equation \eqref{eq:RXHoms} and the identities:
\begin{align*}
R[X](j, \omega) &= f_\omega R[X] f_j \simeq Xe_j = X(j), \\
R[X](\omega, j) &= f_j R[X] f_\omega = 0, 
\end{align*}
for $j = 1, \dots, n$. To prove naturality, given $\alpha: j \to k$, it is clear the following diagram commutes.
\begin{equation*}
\xymatrix{
R[X](k, \omega) \ar[r]^{\;\;\sim} \ar[d]^{(\alpha, -)} & Xe_k \ar[d]^{-\cdot\alpha} \ar@{=}[r] & X(k) \ar[d]^{X(\alpha)} \\
R[X](j, \omega) \ar[r]^{\;\;\sim} & Xe_j \ar@{=}[r] & X(j)
}
\end{equation*}
Checking the remaining bijections are natural is straightforward.
\end{proof}
\end{lemma}

\begin{proposition}\label{prop:OnePointExt}
Let $F: R \to R[X]$ be the inclusion functor of $R$ into a one-point extension $R[X]$. The left adjoint $F_L$  to $\mathrm{res}_F$ exists and is defined pointwise by
\begin{equation*}
F_L M(i) = \begin{cases}
M(i) & \text{if } i \ne \omega,\\
0 & \text{if } i = \omega,
\end{cases}
\end{equation*}
and coincides with the zero embedding $\mathrm{res}_\pi$. The right adjoint $F_R$ to $\mathrm{res}_F$ also exists and is defined pointwise by
\begin{equation*}
F_R M(i) = \begin{cases}
M(i) & \text{if } i \ne \omega,\\
\Hom_R(X, M) & \text{if } i = \omega.
\end{cases}
\end{equation*}
Both adjoints are full and faithful interpretation functors. They coincide on modules $M \in \ModR{R}$ satisfying $\Hom_R(X, M) = 0$ and preserve any Auslander-Reiten sequence in $\modR{R}$ that begins on such a module.
\begin{proof}
The adjoints are constructed by Proposition~\ref{prop:KanExt} (i) (ii) and simplified using Lemma~\ref{lemma:RXProjRes} and the (co-)Yoneda isomorphisms. Moreover, by Lemma~\ref{lemma:RXProjRes}, the restriction of all representable functors are finitely-presented (note this requires that $X = \mathrm{res}_F R[X](-, \omega)$ be finitely-presented) and so, by Proposition~\ref{prop:KanExt} (vi) (vii),  $F_L$ commutes with products, $F_R$ commutes with direct limits, and both are interpretation functors. Further still, since $F$ is an inclusion, $F_L$ and $F_R$ are full and faithful by Proposition~\ref{prop:KanExt} (v).

Note $\pi$ is left adjoint to $F$, indeed, we have already established the bijections $R(\pi(i), j) \simeq R[X](i, F(j))$ which are easily checked to be natural. Therefore, $\pi^{\op}$ is \emph{right} adjoint to $F^{\op}$ (we must be careful and use ``$\op$'' here) and the adjunction $F^{\op} \dashv \pi^{op}$ lifts to an adjunction $\mathrm{res}_{\pi^{\op}} \dashv \mathrm{res}_{F^{\op}}$ by Corollary~\ref{cor:AdjunctLift}. Therefore (after forgoing the ``$\op$''s once again) we have $\mathrm{res}_\pi \simeq F_L$ by uniqueness of adjoints.


The functors $F_L$ and $F_R$ -- restricted to $\modR{R} \to \modR{R[X]}$ -- are just the two embeddings defined in \cite[\S XV.1]{simson_skow1} and the final statement is \cite[XV.1.7]{simson_skow1}.
\end{proof}
\end{proposition}

Similarly, for one-point coextensions, we have a canonical inclusion functor $G: R \to [X]R$, and a projection $\tau: [X]R \to R$ with corresponding zero embedding $\mathrm{res}_\tau: \ModR{R} \to \ModR{[X]R}$.

\begin{proposition}\label{prop:OnePointCoExt}
Let $G: R \to [X]R$ be the inclusion functor of $R$ into a one-point coextension. The left adjoint $G_L$  to $\mathrm{res}_G$ exists and is defined by
\begin{equation*}
G_L M(i) = \begin{cases}
M(i) & \text{if } i \ne \omega,\\
M \otimes_R X^\star & \text{if } i = \omega.
\end{cases}
\end{equation*}
The right adjoint $G_R$ to $\mathrm{res}_G$ exists and is defined by
\begin{equation*}
G_R M(i) = \begin{cases}
M(i) & \text{if } i \ne \omega,\\
0 & \text{if } i = \omega,
\end{cases}
\end{equation*}
and coincides with the zero embedding $\mathrm{res}_\tau$.
Both adjoints are full and faithful interpretation functors. They coincide on modules $M \in \ModR{R}$ satisfying $\Hom_R(M, X) = 0$ and preserve any Auslander-Reiten sequence in $\modR{R}$ that ends on such a module.
\begin{proof}
The is proof similar to Proposition~\ref{prop:OnePointExt} but we point out the main differences. Here $X^\star = \mathrm{res}_G [X]R(\omega, -)$ is finitely-presented, since $X$ is, and this is required for the left adjoint $G_L$ to commute with products. Since $\tau$ is right adjoint to $G$ (i.e. $\tau^{\op}$ is left adjoint to $G^{\op}$) we get $\mathrm{res}_{G^{\op}}$ left adjoint to $\mathrm{res}_{\tau^{\op}}$. Hence (forgoing the ``$\op$'') $\mathrm{res}_\tau \simeq G_R$. Finally, we note that
\begin{equation*}
(M \otimes_R X^\star)^\star = \Hom_\field(M \otimes_R X^\star, \field) \simeq \Hom_R(M, (X^\star)^\star) \simeq \Hom_R(M, X)
\end{equation*}
as $X$ is finitely-presented, so $\Hom_R(M, X) = 0$ if and only if $M \otimes_R X^\star = 0$. The final statement is the dual of \cite[XV.1.7]{simson_skow1}.
\end{proof}
\end{proposition}

We later make use of the fact that $F_R: \ModR{R} \to \ModR{R[X]}$ commutes with direct limits. The functor $G_L: \ModR{R} \to \ModR{[X]R}$ doesn't commute with \emph{all} inverse limits, but we only need the following partial result.

\begin{corollary}\label{prop:CoExtInverseLimits}
If $((M_i)_{i \in \mathbb{N}}, (\alpha_{ji}: M_j \to M_i)_{j \ge i}$ is an inverse system in $\modR{R}$, then $G_L(\varprojlim M_i) \simeq \varprojlim G_L(M_i)$.
\begin{proof}
In the notation of Proposition~\ref{prop:KanExt}, the left adjoint $G_L$ commutes with the inverse limit if and only if for all $v \in [X]R$ the canonical map
\begin{equation*}
\left(\varprojlim M_i\right) \otimes_R [X]R(v, G-) \to \varprojlim \left(M_i \otimes_R [X]R(v, G-)\right)
\end{equation*}
is a $\field$-linear isomorphism \cite[V.4, Ex.5]{maclane}. For $v \ne \omega$, we have $[X]R(v, G-) = R(v, -)$ and, under the co-Yoneda isomorphisms, the above map becomes the obvious identity
\begin{equation*}
\left(\varprojlim M_i\right)(v) \to \varprojlim M_i(v)
\end{equation*}
whereas, for $v = \omega$, we have $[X]R(\omega, G-) = X^\star$ and the above map becomes the canonical map
\begin{equation*}
\left(\varprojlim M_i\right) \otimes_R X^\star \to \varprojlim (M_i \otimes_R X^\star)
\end{equation*}
which is an isomorphism by Proposition~\ref{prop:TensInvs}.
\end{proof}
\end{corollary}

\section{Closure of Some Non-stable Tubes}\label{sec:NonStableTubes}
In this section we will construct a class of algebras containing a separating tubular family with a single non-stable tube. We will compute the Ziegler closure of this tube, describing the infinite-dimensional points therein in terms of rays and corays.

For the precise definition of a tube we refer to \cite[\S 1]{deste_ringel}. Every arrow in a tube either ``points to infinity'' or ``points to the mouth''. A \textbf{ray} is an infinite path $X[1] \to X[2] \to X[3] \to \cdots$ in a tube, with pairwise distinct vertices and all arrows pointing to infinity. A vertex $X$ in a tube is a \textbf{ray vertex} if there exists an infinite \textbf{sectional path} $X[1] \to X[2] \to \cdots$ (i.e. $X[i] \ne \tau X[i+2]$ for all $i \ge 1$) starting at $X = X[1]$ and containing every sectional path starting at $X$ (such an infinite path $X[1] \to X[2] \to \cdots$ is indeed a ray \cite[4.6.3]{ringel}). In this terminology, not all rays need start on a ray vertex (e.g. in Example~\ref{ex:Gamma232} below, $X_0[1] \to X_0[2] \to \cdots $ is ray but $X_0[1]$ is not a ray vertex). A \textbf{coray} and \textbf{coray vertex} are defined dually.

By a \textbf{ray insertion} at a ray vertex $X$, we mean a 1-fold ray insertion at $X$ --- the definition of a general \textbf{$n$-fold ray insertion} is \cite[2.1]{deste_ringel}. For example, a single ray insertion at ray vertex $Y[1]$ of the quiver below left, results in a new ray $Y^\prime[1] \to Y^\prime[2] \to \cdots$ inserted as depicted in the quiver below right. 

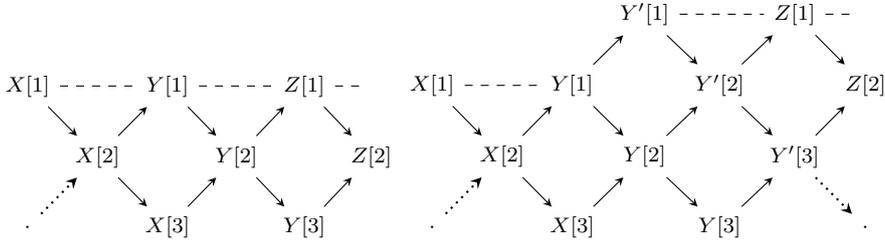
\begin{equation*}
\begin{tikzpicture}[every node/.style={anchor=center, font=\scriptsize},>=stealth]
\matrix [matrix of math nodes, row sep=0.4cm, column sep=0.08cm, anchor=south east] (a) at (0,0) {
X[1] && Y[1] && Z[1] & \, \\
& X[2] && Y[2] && Z[2] \\
\cdot && X[3] && Y[3] \\
};
\matrix [matrix of math nodes, row sep=0.4cm, column sep=0.08cm, anchor=south east] (b) at (6.5,0) {
&&& Y^\prime[1] && Z[1] & \, \\
X[1] && Y[1] && Y^\prime[2] && Z[2] \\
& X[2] && Y[2] && Y^\prime[3] \\
\cdot && X[3] && Y[3] && \cdot  \\
};
\draw[dashed] (a-1-1) to (a-1-3);
\draw[dashed] (a-1-3) to (a-1-5);
\draw[dashed] (a-1-5) to (a-1-6);
\draw[->] (a-1-1) to (a-2-2);
\draw[->] (a-2-2) to (a-1-3);
\draw[->] (a-1-3) to (a-2-4);
\draw[->] (a-2-4) to (a-1-5);
\draw[->] (a-1-5) to (a-2-6);
\draw[thick, dotted, ->] (a-3-1) to (a-2-2);
\draw[->] (a-2-2) to (a-3-3);
\draw[->] (a-3-3) to (a-2-4);
\draw[->] (a-2-4) to (a-3-5);
\draw[->] (a-3-5) to (a-2-6);

\draw[dashed] (b-1-4) to (b-1-6);
\draw[dashed] (b-1-6) to (b-1-7);
\draw[dashed] (b-2-1) to (b-2-3);
\draw[->] (b-2-3) to (b-1-4);
\draw[->] (b-1-4) to (b-2-5);
\draw[->] (b-2-5) to (b-1-6);
\draw[->] (b-1-6) to (b-2-7);
\draw[->] (b-2-1) to (b-3-2);
\draw[->] (b-3-2) to (b-2-3);
\draw[->] (b-2-3) to (b-3-4);
\draw[->] (b-3-4) to (b-2-5);
\draw[->] (b-2-5) to (b-3-6);
\draw[->] (b-3-6) to (b-2-7);
\draw[thick, dotted, ->] (b-4-1) to (b-3-2);
\draw[->] (b-3-2) to (b-4-3);
\draw[->] (b-4-3) to (b-3-4);
\draw[->] (b-3-4) to (b-4-5);
\draw[->] (b-4-5) to (b-3-6);
\draw[thick, dotted, ->] (b-3-6) to (b-4-7);

\end{tikzpicture}
\end{equation*}

\paragraph{}
Let us first define the translation quivers for the type of tube we will deal with. For integers $p \ge 0$ and $0 \le m \le n$ define $\mathbf{\Gamma(p,n,m)}$ to be the translation quiver with the following vertices for all $j \in \mathbb{N}$ (interpret an empty interval as having no vertices of that kind).
\begin{align*}
X_i[j] & \text{ for } 0 \le i \le n-m,\\
Y_i[j] & \text{ for } 1 \le i \le m,\\
Z_i[j] & \text{ for } 1 \le i \le p.
\end{align*}
The arrows of $\Gamma(p,n,m)$ are as follows, for all $j \in \mathbb{N}$.
\begin{align*}
\text{for } 0 \le i \le n-m &\quad X_i[j] \to X_i[j+1],\\
\text{for } 1 \le i \le n-m &\quad X_{i-1}[j] \to X_i[j],\\
\text{if } m=p=0 &\quad X_n[j+1] \to X_0[j].
\end{align*}
Also if $m \ge 1$, then
\begin{align*}
&X_{n-m}[j] \to Y_1[j], \\
\text{ for } 1 \le i \le m \quad& Y_i[j] \to Y_i[j+1], \\
\text{ for }  1 \le i < m \quad& Y_i[j+1] \to Y_{i+1}[j], \\
\text{ if } p \ge 1 \quad& Y_m[j+2] \to Z_1[j], \\
\text{ if } p = 0 \quad& Y_m[j+2] \to X_0[j].
\end{align*}
Also if $p \ge 1$, then
\begin{align*}
&Z_p[j+1] \to X_0[j] \\
\text{ for } 1 \le i \le p \quad& Z_i[j] \to Z_i[j+1], \\
\text{ for } 1 \le i < p \quad& Z_i[j+1] \to Z_{i+1}[j], \\
\text{ if } m=0 \quad& X_n[j+1] \to Z_1[j].
\end{align*}
The translation of $\Gamma(p,n,m)$ is as follows, for all $j \in \mathbb{N}$. 
\begin{equation*}
\tau X_i[j+1] = X_{i-1}[j] \text{ for } 1 < i \le m-n
\end{equation*}
and
\begin{equation*}
\tau X_0[j] = \begin{cases}
X_n[j] &\text{if } m=0,\,p=0,\\
Y_n[j+1] &\text{if } m\ge 1,\,p \ge 1,\\
Z_m[j] &\text{if } m \ge 1,\,p=0.
\end{cases}
\end{equation*}
Also if $m \ge 1$, then
\begin{equation*}
\tau Y_i[j+1] = \begin{cases}
X_{n-m}[j] & \text{if } i=1,\\
Y_{i-1}[j+1] & \text{if } 1 < i \le m.
\end{cases}
\end{equation*}
Also if $p \ge 1$, then $\tau Z_i[j] = Z_{i-1}[j] \text{ for } 1 < i \le p$ and
\begin{equation*}
\tau Z_1[j] = \begin{cases}
Y_m[j+1] &\text{ if } m \ge 1,\\
X_n[j] &\text{ if } m=0.
\end{cases}
\end{equation*}

\begin{example}\label{ex:Gamma232}
The translation quiver $\Gamma(2,3,2)$ is depicted (near the mouth) in the following diagram. Only the translates at the mouth are shown, by a dashed line (note $Y_1[1]$ and $Y_2[1]$ are projective-injective); the left and right edges are to be identified.
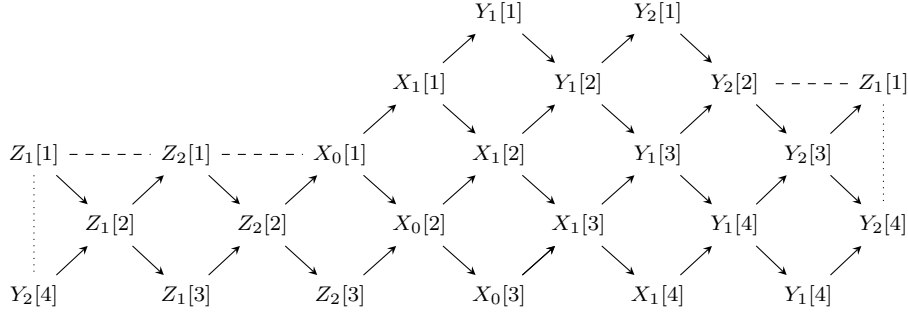
\begin{equation*}
\begin{tikzpicture}[every node/.style={anchor=center, font=\scriptsize},>=stealth]
\matrix [matrix of math nodes, row sep=0.4cm, column sep=0.08cm, anchor=south east] (a) at (0,0) {
&&&&&& Y_1[1] && Y_2[1] \\
&&&&& X_1[1] && Y_1[2] && Y_2[2] && Z_1[1] \\
Z_1[1] && Z_2[1] && X_0[1] && X_1[2] && Y_1[3] && Y_2[3] \\
& Z_1[2] && Z_2[2] && X_0[2] && X_1[3] && Y_1[4] && Y_2[4] \\
Y_2[4] && Z_1[3] && Z_2[3] && X_0[3] && X_1[4] && Y_1[4]  \\
};
\draw[->] (a-3-1) to (a-4-2);
\draw[->] (a-4-2) to (a-3-3);
\draw[->] (a-3-3) to (a-4-4);
\draw[->] (a-4-4) to (a-3-5);
\draw[->] (a-3-5) to (a-4-6);
\draw[dashed] (a-3-1) to (a-3-3);
\draw[dashed] (a-3-3) to (a-3-5);
\draw[->] (a-4-6) to (a-3-7);
\draw[->] (a-3-7) to (a-4-8);
\draw[->] (a-4-8) to (a-3-9);
\draw[->] (a-3-9) to (a-4-10);
\draw[->] (a-4-10) to (a-3-11);
\draw[->] (a-3-5) to (a-2-6);
\draw[->] (a-2-6) to (a-1-7);
\draw[->] (a-1-7) to (a-2-8);
\draw[->] (a-2-6) to (a-3-7);
\draw[->] (a-3-7) to (a-2-8);
\draw[->] (a-2-8) to (a-3-9);
\draw[->] (a-2-8) to (a-1-9);
\draw[->] (a-3-9) to (a-2-10);
\draw[dashed] (a-2-12) to (a-2-10);
\draw[->] (a-1-9) to (a-2-10);
\draw[->] (a-2-10) to (a-3-11);
\draw[->] (a-3-11) to (a-2-12);
\draw[->] (a-3-11) to (a-4-12);
\draw[->] (a-5-1) to (a-4-2);
\draw[->] (a-4-2) to (a-5-3);
\draw[->] (a-5-3) to (a-4-4);
\draw[->] (a-4-4) to (a-5-5);
\draw[->] (a-5-5) to (a-4-6);
\draw[->] (a-4-6) to (a-5-7);
\draw[->] (a-5-7) to (a-4-8);
\draw[->] (a-5-7) to (a-4-8);
\draw[->] (a-4-8) to (a-5-9);
\draw[->] (a-5-9) to (a-4-10);
\draw[->] (a-4-10) to (a-5-11);
\draw[->] (a-5-11) to (a-4-12);
\draw[dotted] (a-3-1) to (a-5-1);
\draw[dotted] (a-2-12) to (a-4-12);
\end{tikzpicture}
\end{equation*}
\end{example}

Further examples are at the end of this section.

\begin{lemma}\label{lemma:G-pnm}
The translation quiver $\Gamma(p,n,m)$ is a tube and is obtained from a stable tube of rank $p+1$ by $n$ ray insertions followed by $m$ coray insertions, in the following way.
\begin{enumerate}[(i)]
\item The quiver $\Gamma(p,0,0)$ is a stable tube of rank $p+1$.
\item $X_n[1]$ is a ray vertex in $\Gamma(p,n,0)$.
\item If $n \ge 1$, then $\Gamma(p,n,0)$ is obtained from $\Gamma(p,n-1,0)$ by a ray insertion at $X_{n-1}[1]$.
\item $X_{n-m}[1]$ is a coray vertex in $\Gamma(p,n,m)$.
\item If $m \ge 1$, then $\Gamma(p,n,m)$ is obtained from $\Gamma(p,n,m-1)$ by a coray insertion at $X_{n-m+1}[1]$.
\end{enumerate}
\begin{proof}

We only note that the $n$ consecutive 1-fold ray insertions, used to obtained $\Gamma(p,n,0)$ from $\Gamma(p,0,0)$, do not equate to a single $n$-fold ray insertion (as defined in \cite[2.1]{deste_ringel}).
\end{proof}
\end{lemma}

\paragraph{Remarks}
\begin{enumerate}[(i)]
\item The vertices of $\Gamma(p,n,m)$ can be partitioned by $p+n+1$ (maximal) rays beginning at the vertices $X_0[1], \dots, X_{n-m}[1]$, $Y_1[1], \dots, Y_m[1]$, and $Z_1[1], \dots, Z_p[1]$. Also, the vertices can be partitioned by $p+m+1$ (maximal) corays ending at the vertices $Y_1[1], \dots, Y_m[1],\, Y_m[2]$ and $Z_1[1], \ldots, Z_p[1]$ if $m \ge 1$, otherwise ending at $X_n[1]$ and $Z_1[1], \dots, Z_p[1]$.
\item The vertices $X_1[1], \dots, X_{n-m}[1]$ are all projective and the vertices $Y_1[1],\dots,Y_m[1]$ are both projective and injective (when they exist); all other vertices are neither projective nor injective.
\item If $n \ge 1$ and $m=n$ the quiver $\Gamma(p,n,n)$ has $n$ projective-injective vertices $Y_1[1],\dots,Y_n[1]$ and becomes a stable tube of rank $p+n+1$ upon removing these points. 
\end{enumerate}

We now define algebras containing tubes of the form $\Gamma(p,n,m)$.

Let $A$ be a tame hereditary algebra with $S$ a simple regular $A$-module. Define the algebra $A[S,n]$ for $n \ge 0$ by iterated one-point extensions as follows.
\begin{equation*}
A[S,n] := \begin{cases}
A & \text{if $n=0$,} \\
(A[S,n-1])[P_{n-1}] & \text{if $n \ge 1$},
\end{cases}
\end{equation*}
where $P_0 = S$ and for $i \ge1$, $P_i = A[S,i](-, \omega_i)$ is the indecomposable projective $A[S,i]$-module corresponding to the extension vertex $\omega_i$ of the one-point extension from $A[S,i-1]$.

Now define the algebra $A[S,n,m]$ for $0 \le m \le n$ by iterated one-point coextensions as follows.
\begin{equation*}
A[S,n,m] := \begin{cases}
A[S,n] & \text{if $m = 0$,} \\
[W_{m-1}](A[S,n,m-1]) & \text{if $m \ge 1$},
\end{cases}
\end{equation*}
where $W_{m-1} := \tau^{m-1} P_n$ is the $(m-1)^{\text{th}}$ translate of $P_n$ in $\modR{A[S,n,m-1]}$ (that such a module exists is proven below).

\paragraph{}
Let $R$ be an algebra and $\mcal{C}$ a tube of $\modR{R}$. An indecomposable module $M \in \mcal{C}$ is a \textbf{ray module} if the corresponding point of $\Gamma(\mcal{C})$ is a ray vertex. One-point extensions by ray modules (and dually, one-point coextensions by coray modules) in standard tubes are a special case of \cite[4.5.1]{ringel}, which used inductively gives the following proposition. Nevertheless, we try to give sufficient detail of the procedure to help clarify subsequent results.

%

\begin{proposition}\label{prop:AStube}
Let $A$ be a tame hereditary algebra with $S$ a simple regular $A$-module. Given integers $0 \le m \le n$ we have
\begin{equation*}
\modR{A[S,n,m]} = \mcal{P} \vee \mcal{T} \vee \mcal{Q}
\end{equation*}
with $\mcal{T}$ a standard $\mathbb{P}^1(\field)$-tubular family separating $\mcal{P}$ from $\mcal{Q}$. The module $S$ lies in a tube $\mathbb{T}$ of the form $\Gamma(\mathbb{T}) = \Gamma(t-1,n,m)$ in $\mcal{T}$, where $t$ is the $\tau_A$-periodicity of $S$. The remaining tubes of $\mcal{T}$ consist of regular $A$-modules and are just the tubes of $\modR{A}$ not containing $S$.
\begin{proof}
We proceed by induction on $0 \le m \le n$. For $m = n = 0$ we have $A[S,0,0] = A$ is a tame hereditary algebra. So $\mcal{T}$ is the class of regular $A$-modules and, by choice of $S$, we know $\mathbb{T}$ is a stable tube of rank $t$. Hence, we may identify $\Gamma(\mathbb{T}) = \Gamma(t-1,0,0)$ by Lemma~\ref{lemma:G-pnm}, and do so in such a way that $S$ corresponds to the vertex $X_0[1]$.

For $m=0$ we proceed inductively on $n \ge 1$. By definition $A[S,n,0] = A[S,n]$ is the one-point extension of $A[S,n-1]$ at $P_{n-1}$. Let $\mathbb{S}$ be the tube of $A[S,n-1]$ containing $S$. Assume that we can identify $\Gamma(\mathbb{S}) = \Gamma(t-1, n-1, 0)$ in such a way that $P_{n-1}$ corresponds to $X_{n-1}[1]$ --- this holds when $n=1$, since $P_0 = S$. Then, by Lemma~\ref{lemma:G-pnm}, $P_{n-1}$ is a ray module and, by \cite[4.5.1]{ringel}, the component of $\modR{A[S,n]}$ containing $P_{n-1}$ is a tube $\mathbb{T}$ obtained from $\mathbb{S}$ by a single ray insertion: $\Gamma(\mathbb{T})$ is obtained from $\Gamma(\mathbb{S})$ by a ray insertion at $X_{n-1}[1]$ and we may identify $\Gamma(\mathbb{T}) = \Gamma(t-1,n,0)$ in such a way that $P_n$ corresponds to $X_n[1]$ (the only ``new'' projective vertex). The module $S$ lies in $\mathbb{T}$ and the remaining claims follow from \cite[4.5.1]{ringel}.

For $n \ge 1$ proceed inductively on $1 \le m \le n$. By definition $A[S,n,m]$ is the one-point coextension of $A[S,n,m-1]$ at $W_{m-1} = \tau^{m-1}P_n$. Let $\mathbb{S}$ be the tube of $A[S,n,m-1]$ containing $S$. Assume that $W_{m-1}$ exists and that we can identify $\Gamma(\mathbb{S}) = \Gamma(t-1,n,m-1)$ in such a way that $W_{m-1}$ corresponds to $X_{n-m+1}[1]$ --- this holds when $m=1$, as $W_0 = \tau^0 P_n = P_n$. Then, by Lemma~\ref{lemma:G-pnm}, $W_{m-1}$ is a coray module and by (the dual of) \cite[4.5.1]{ringel}, the component of $\modR{A[S,n,m]}$ containing $W_{m-1}$ is a tube $\mathbb{T}$ obtained from $\mathbb{S}$ by a single coray insertion: $\Gamma(\mathbb{T})$ is obtained from $\Gamma(\mathbb{S})$ by a coray insertion at $X_{n-m+1}[1]$ and we may identify $\Gamma(\mathbb{T}) = \Gamma(t-1,n,m)$ in such a way that $W_{m-1}$ corresponds to $Y_1[2]$, so that $W_m = \tau^m P_n= \tau W_{m-1}$ exists and corresponds to $\tau Y_1[2] = X_{n-m}[1]$. The module $S$ lies in $\mathbb{T}$ and the remaining claims follow from \cite[4.5.1]{ringel}.
\end{proof}
\end{proposition}

\begin{lemma}\label{lemma:RayZero}
If $M$ is a ray module in a standard tube $\mathbb{T}$, then for all indecomposable $X \in \mathbb{T}$ we have $(M, X) = 0$ unless $X \simeq M[i]$ for some $i \ge 1$.
\begin{proof}
This is contained in the proof of \cite[4.5.1]{ringel}. There, it is shown that any non-sectional path beginning at $M$ is equal to zero (modulo the mesh relations).
\end{proof}
\end{lemma}

Let $R$ be an algebra and $\mathbb{T}$ a tube of $\modR{R}$. Given a chain of irreducible morphisms $U[1] \to U[2] \to \cdots$ representing a maximal ray in $\Gamma(\mathbb{T})$, the direct limit $U[\infty] := \varinjlim_j U[j]$ is called the \textbf{$U$-pr\"ufer} module. Similarly, given a chain of irreducible morphisms $\cdots \to [2]U \to [1]U$ representing a maximal coray in $\Gamma(\mathbb{T})$, the inverse limit $[\infty]U := \varprojlim_j U[j]$ is called the \textbf{$U$-adic} module. Using the standard duality $(-)^\star: (\modR{A})^{\op} \to \modR{A^{\op}}$, a ray (resp. coray) of $A$-modules is sent to a coray (resp. ray) of $A^{\op}$-modules. Furthermore
\begin{equation*}
(U[\infty])^\star = \Hom_\field(\varinjlim U[i], \field) \simeq \varprojlim \Hom_\field(U[i], \field) = \varprojlim U[i]^\star
\end{equation*}
so the $\field$-dual of a pr\"ufer module is an adic module, and every adic module is of this form.

For instance, if $R = A[S,n,m]$ and $\mathbb{T}$ is the tube containing $S$ as in the preceding proposition, then by the remarks following Lemma~\ref{lemma:G-pnm}, to $\mathbb{T}$ we can associate $t+n$ pr\"ufer modules and $t+m$ adic modules (to do so unambiguously requires showing these modules are independent of the choice of irreducible maps; this is a consequence of the following theorem). Proposition~\ref{prop:AStube} describes the tube $\mathbb{T}$ and we now describe the closure $\mathrm{cl}(\mathbb{T})$ of $\mathbb{T}$ in the Ziegler spectrum $\Zg_R$. 


\begin{lemma}\label{lemma:RayZgZero}
Let $A$ be a tame hereditary algebra, with $\mathbb{T}$ a tube in $\modR{A}$ and $M \in \mathrm{cl}(\mathbb{T})$ in the Ziegler closure of $\mathbb{T}$.
\begin{enumerate}[(i)]
\item If $X \in \mathbb{T}$ is a ray module, then $(X, M) = 0$ unless $M \simeq X[j]$ for some $j \in \mathbb{N} \cup \{ \infty \}$. 
\item If $X \in \mathbb{T}$ is a coray module, then $(M, X) = 0$ unless $M \simeq [j]X$ for some $j \in \mathbb{N} \cup \{ \infty \}$.
\end{enumerate}
\begin{proof}
For $M$ finite-dimensional the result follows from Lemma~\ref{lemma:RayZero} (and its dual), so suppose $M$ is infinite-dimensional. Then, by \cite{ringel2}, $M$ is one of the pr\"ufer or adic modules of $\mathbb{T}$ or the generic module.

Suppose $X$ is a ray module. By \cite[X.2.8]{simson_skow2} a non-zero map $X \to Y$, with $Y \in \indR{T}$ but $Y$ not in $\mathbb{T}$, factors through the composition of morphisms $X = X[1] \to X[2] \to \cdots \to X[n]$, from the ray beginning at $X$, for arbitrarily large $n \ge 1$. By \cite[5.3.31]{prest} the same result holds true in case $Y = M$. Thus, a non-zero map $f: X \to M$ induces a non-zero map $X[\infty] \to M$. However, it is shown in \cite{reiten_ringel} that $(X[\infty], N) = 0$ whenever $N \in \mathbb{T}$ or $N = G$ (the generic module), and consequently
\begin{equation*}
(X[\infty], [\infty]V) = (X[\infty], \varprojlim [j]V) = \varprojlim (X[\infty], [j]V) = 0
\end{equation*}
for any adic module $[\infty]V$. Thus, if $(X, M) \ne 0$, then $M$ must be a pr\"ufer module. If $M = U[\infty]$ for $U \ne X$, then as $X$ is finitely-presented we have
\begin{equation*}
(X, M) = (X, \varinjlim U[j]) = \varinjlim (X, U[j]) = 0
\end{equation*}
this completes the proof of (i).

Suppose $X$ is a coray module, then $(M, X) = 0$ for any pr\"ufer or generic module $M$ by \cite{reiten_ringel} again. So suppose $M = [\infty]V$  is an adic module with $V \ne X$. We have $(M, X) \simeq (M \otimes_R X^\star)^\star$ as $X$ is finite-dimensional, and by Proposition~\ref{prop:TensInvs} we know
\begin{equation*}
M \otimes_R X^\star = \left(\varprojlim [i]V \right) \otimes_R X^\star \simeq \varprojlim \left( [i]V \otimes_R X^\star \right) \simeq \varprojlim\left( ([i]V, X)^\star \right) = 0
\end{equation*}
since $([i]V, X) \simeq (X^\star, [i]V^\star) = 0$ for all $i \ge 1$, by standard duality and the fact that $[i]V^\star$ belongs to a ray not beginning at $X^\star$. Thus $(M, X) = 0$ as claimed.
\end{proof} 
\end{lemma}

\begin{theorem}\label{prop:ASTubeClosure}
Let $A$ be a tame hereditary algebra with $S$ a simple regular $A$-module of $\tau_A$-periodicity $t$. Given integers $0 \le m \le n$ let $\mathbb{T}$ be the tube of $\modR{A[S,n,m]}$ containing $S$, as in Proposition~\ref{prop:AStube}. Then the $t+n$ pr\"ufer and $t+m$ adic modules of $\mathbb{T}$ are indecomposable, pure-injective, and pairwise non-isomorphic. These, together with a unique generic module, comprise all infinite-dimensional points in the Ziegler closure $\mathrm{cl}(\mathbb{T})$ of $\mathbb{T}$ in $\Zg_{A[S,n,m]}$.
\begin{proof}
The proof proceeds inductively as in Proposition~\ref{prop:AStube}.

If $n = m = 0$, then $A[S,0,0] = A$ is a tame hereditary algebra. Identify $\Gamma(\mathbb{T}) = \Gamma(t-1,0,0)$ such that $S$ represents $X_0[1]$. The $X_0$-pr\"ufer and $X_0$-adic are just the $S$-pr\"ufer and $S$-adic of \cite{ringel2} (cf. \cite{prest2}). Similarly, the $Z_i$-pr\"ufer and $Z_i$-adic are the $\tau^{t-i}S$-pr\"ufer and $\tau^{t-i} S$-adic respectively, for $i = 1, \dots, t-1$. The remaining claims follow from \cite{ringel2}.

For $m = 0$ we proceed inductively on $n \ge 1$. As defined in Proposition~\ref{prop:OnePointExt}, the one-point extension $A[S,n] = (A[S,n-1])[P_{n-1}]$ comes with two embeddings $\ModR{A[S,n-1]} \to \ModR{A[S,n]}$  --- namely, let $F_L$ (resp. $F_R$) denote the left (resp. right) adjoint to the restriction functor $\ModR{A[S,n]} \to \ModR{A[S,n-1]}$. Let $\mathbb{S}$ be the tube of $\modR{A[S,n-1]}$ containing $S$ and $\mathbb{T}$ be the tube of $\modR{A[S,n]}$ containing $S$. Then $F_L$ (the zero embedding) restricts to a map $\mathbb{S} \to \mathbb{T}$ compatible with identifying the points of $\Gamma(t-1, n-1, 0) = \Gamma(\mathbb{S})$ with the correspondingly named points of $\Gamma(t-1,n ,0) = \Gamma(\mathbb{T})$. Now $F_R$ coincides with $F_L$ for all indecomposable $M \in \mathbb{S}$ satisfying $(X_{n-1}[1], M) = 0$ by Proposition~\ref{prop:OnePointExt}. This is the case unless $M \simeq X_{n-1}[j]$ for some $j \ge 1$, by Lemma~\ref{lemma:RayZero}. Thus $F_L$ and $F_R$ differ only on the ray beginning at $X_{n-1}[1]$, where $F_R X_{n-1}[j] = X_n[j]$ for $j \ge 1$ (consult the proof of \cite[4.5.1]{ringel}). In this way, $F_R$ also restricts to a map $\mathbb{S} \to \mathbb{T}$.

As $F_L$ and $F_R$ are full and faithful interpretation functors, they induce closed embeddings $\Zg_{A[S,n-1]} \to \Zg_{A[S,n]}$ by Proposition~\ref{prop:InterpFunc}. Let $\mathrm{cl}(\mathbb{S})$ denote the closure of $\mathrm{ind}(\mathbb{S})$ in $\Zg_{A[S,n-1]}$ and set $C_1 = F_L(\mathrm{cl}(\mathbb{S}))$ and $C_2 = F_R(\mathrm{cl}(\mathbb{S}))$. Since both $F_L$ and $F_R$ are closed maps, $C_1 \cup C_2$ is a closed subset of $\Zg_{A[S,n]}$ containing all points of $\mathrm{ind}(\mathbb{T})$, hence $\mathrm{cl}(\mathbb{T}) \subseteq C_1 \cup C_2$. However, as $F_L(\mathbb{S}) \subseteq \mathbb{T}$, by continuity we have $C_1 = F_L(\mathrm{cl}(\mathbb{S})) \subseteq \mathrm{cl}(F_L(\mathbb{S})) \subseteq \mathrm{cl}(\mathbb{T})$ and similarly for $C_2$. Therefore $\mathrm{cl}(\mathbb{T}) = C_1 \cup C_2$.

Now choose irreducible morphisms in $\mathbb{S}$ satisfying the mesh relations of $\Gamma(\mathbb{S})$ (recall $\mathbb{S}$ is standard). These morphisms remain irreducible in $\mathbb{T}$ except for those of the form $X_{n-1}[j+1] \to Z_1[j]$ which factor as a composition of irreducible morphisms $X_{n-1}[j+1] \to X_n[j+1] \to Z_1[j]$. These maps, together with the maps $X_n[j] \to X_n[j+1] = F_R(X_{n-1}[j] \to X_{n-1}[j])$ and one additional map $X_{n-1}[1] \to X_n[1]$ (the inclusion of $P_{n-1} = \mathrm{rad}(P_n)$ into $P_n$), form a system of irreducible morphisms satisfying the mesh relations of $\Gamma(\mathbb{T})$ (again, consult the proof \cite[4.5.1, p.219]{ringel}). If $t=1$, then read ``$Z_i$'' as ``$X_0$'' in this and the following paragraphs.

The $t+n$ maximal rays in $\mathbb{T}$, except the ray beginning at $X_n[1]$, are also rays in $\mathbb{S}$. Thus, the $X_i$-pr\"ufers (for $i=0,\dots,n-1$) and $Z_i$-pr\"ufers (for $i=1,\dots,t-1$) are indecomposable and pure-injective by assumption. For the remaining ray $X_n[1] \to X_n[2] \to \cdots$ we have
\begin{equation*}
X_n[\infty] = \varinjlim X_n[j] = \varinjlim F_R(X_{n-1}[j]) = F_R(\varinjlim X_{n-1}[j]) = F_R(X_{n-1}[\infty])
\end{equation*}
since $F_R$ commutes with direct limits. Therefore the $X_n$-pr\"ufer is also indecomposable and pure-injective, since $F_R$ is a full and faithful interpretation functor.

Now if $\cdots \to [2]V \to [1]V$ is one of the $t$ maximal corays of $\mathbb{T}$, there will be some $k \ge 1$ and $j \ge 2$ such that we get the configuration
\begin{equation*}
\xymatrix{
[k+2]V \ar[r] \ar@{=}[d] & [k+1]V \ar[r] \ar@{=}[d] & [k]V \ar@{=}[d] \\
X_{n-1}[j] \ar[r] & X_n[j] \ar[r] & Z_1[j-1]
}
\end{equation*}
In fact, this will happen infinitely often and periodically (also if $V = X_n$, then $[2]V \to [1]V = X_{n-1}[1] \to X_n[1]$). But as noted above, the composition $X_{n-1} \to Z_1[j-1]$ is an irreducible morphism in $\mathbb{S}$. The inverse limit along the coray is not changed if we remove all points of the form $X_n[j]$ for $j \ge 1$. This results in a coray of $\mathbb{S}$ and so, since $F_R$ preserves indecomposable pure-injectives, the $V$-adic $[\infty]V = \varprojlim [j]V$ is again indecomposable and pure-injective.

Together with the generic module (which remains generic over $A[S,n,m]$), we have described all infinite-dimensional points in $C_1$. The $X_n$-pr\"ufer lies in $C_2$ and we claim $C_2$ contains no points we haven't already considered. If we add to our induction hypothesis the following property (which holds for $n=1$ by Lemma~\ref{lemma:RayZgZero})
\begin{equation}\label{eq:InductiveProperty}
(X_{n-1}[1], M) = 0 \text{ for } M \in \mathrm{cl}(\mathbb{S}) \text{ unless } M \simeq X_{n-1}[j] \text{ for } j \in \mathbb{N} \cup \{ \infty \}
\end{equation}
then our claim follows immediately, for then $F_L$ and $F_R$ coincide on all infinite-dimensional points in $\mathrm{cl}(\mathbb{S})$ except the $X_{n-1}$-pr\"ufer module. To complete the induction we must now show the following.
\begin{equation}\label{eq:InductiveProperty2}
(X_n[1], M) = 0 \text{ for } M \in \mathrm{cl}(\mathbb{T}) \text{ unless } M \simeq X_n[j] \text{ for } j \in \mathbb{N} \cup \{ \infty \}
\end{equation}
So, suppose $M = F_RM^\prime \in C_2$, then
\begin{equation*}
(X_n[1], M) = (F_RX_{n-1}[1], F_RM^\prime) \simeq (X_{n-1}[1], M^\prime)
\end{equation*}
so $(X_n[1], M) = 0$ unless $M^\prime \simeq X_{n-1}[j]$ for some $j \in \mathbb{N} \cup \{ \infty \}$ by \eqref{eq:InductiveProperty}, that is, unless $M = F_RM^\prime \simeq F_RX_{n-1}[j] = X_n[j]$. The only other option is $M = X_{n-1}[j]$ for some $j \in \mathbb{N} \cup \{\infty\}$, but we know $(X_n[1], X_{n-1}[j]) = 0$ for all $j \in \mathbb{N}$ by Lemma~\ref{lemma:RayZero}, hence $(X_n[1], X_{n-1}[\infty]) = 0$, completing the proof of \eqref{eq:InductiveProperty2}.

For $n \ge 1$ we proceed inductively on $1 \le m \le n$. The proof is essentially dual to the case above. We again have two embeddings $F_L, F_R: \ModR{A[S,n,m-1]} \to \ModR{A[S,n,m]}$, given by Proposition~\ref{prop:OnePointCoExt}, although their roles are reversed (here $F_R$ is the zero embedding). Recall $A[S,n,m]$ is obtained from $A[S,n,m-1]$ by a coray insertion. $F_R$ preserves the adic, pr\"ufer and generic modules of the tube containing $S$ in $\modR{A[S,n,m-1]}$. There is one additional adic $A[S,n,m]$-module in the image of $F_L$ (corresponding to the new coray) and here Corollary~\ref{prop:CoExtInverseLimits} is necessary to ensure $F_L$ commutes with the relevant inverse limit.
\end{proof}
\end{theorem}

We extract the following result implicit in the proof of the above theorem.

\begin{proposition}\label{prop:ASTubeTop}
Let $A$ be a tame hereditary algebra with $S$ a simple regular $A$-module. Given integers $0 \le m \le n$, let $\mathbb{T}$ be the tube of $\modR{A[S,n,m]}$ containing $S$, then a subset $\mathfrak{C} \subseteq \mathrm{cl}(\mathbb{T})$ is closed in $\Zg_{A[S,n,m]}$ if and only if the following hold:
\begin{enumerate}[(i)]
\item The $U$-pr\"ufer belongs to $\mathfrak{C}$ whenever $\mathfrak{C}$ contains infinitely many modules along the ray beginning at $U$.
\item The $V$-adic belongs to $\mathfrak{C}$ whenever $\mathfrak{C}$ contains infinitely many modules along the coray ending at $V$.
\item The generic module belongs to $\mathfrak{C}$ whenever $\mathfrak{C}$ contains infinitely-many finite-dimensional points or at least one infinite-dimensional point.
\end{enumerate}
The CB ranks of the points in $\mathrm{cl}(\mathbb{T})$ are as follows:
\begin{enumerate}[(a)]
\item Each finite-dimensional point has CB rank 0.
\item Each pr\"ufer and adic module has CB rank 1.
\item The generic module has CB rank 2.
\end{enumerate}
and furthermore $\mathrm{cl}(\mathbb{T})$ satisfies the isolation condition.
\begin{proof}
The statements (i)---(iii) follow by the induction of Theorem~\ref{prop:ASTubeClosure}: at each step the tube is covered by two tubes whose closure has the stated topology (or more accurately, by the same tube in two different ways). The base case at $n=m=0$ is given by the main theorem of \cite{ringel2}.

Now let $\mathfrak{C} = \mathrm{cl}(\mathbb{T})$. We show each point $M \in \mathfrak{C}$ is isolated in its own closure by an $M$-minimal pp-pair, then the isolation condition holds for $\mathfrak{C}$ by \cite[5.3.16]{prest}. 

If $M \in \mathrm{cl}(\mathbb{T})$ is finite-dimensional, then $\{ M \}$ is both open and closed in $\mathfrak{C}$ by Proposition~\ref{prop:ZgFdAlg}, so $M$ has CB rank 0. Furthermore, $\{ M \}$ is given by a minimal pp-pair by \cite[5.3.31]{prest}.

If $M = U[\infty]$ is a pr\"ufer module, $\{ M \}$ cannot be open in $\mathfrak{C}$ (for else $\mathfrak{C} \backslash \{ M \}$ is a closed set containing all points $U[j]$ for $j \ge 1$ but not $U[\infty]$, contradicting (i)), so $M$ has CB rank $> 0$.  Similarly, the adic and generic modules cannot be isolated in $\mathfrak{C}$, thus $\mathfrak{C}^{(1)} = \mathfrak{C} \backslash \mathbb{T}$ is the set of infinite-dimensional points in $\mathfrak{C}$. Now $\{ M \}$ is open in $\mathfrak{C}^{(1)}$ and so $M$ has CB rank 1. Similar reasoning tells us each adic module has CB rank 1 and $\mathfrak{C}^{(2)} = \{ G \}$, giving the generic module CB rank 2.

At some stage of the induction of Theorem~\ref{prop:ASTubeClosure}, when the pr\"ufer module $M = U[\infty]$ first appears, the equation~\eqref{eq:InductiveProperty2} tells us $(U[1], -)$ isolates $M$ amongst the infinite-dimensional points (at that stage). Although $(U[1], -)$ may not isolate $M$ in $\mathfrak{C}^{(1)}$ (for $U[1]$ embeds into any pr\"ufer module appearing at a later stage), it does isolate $M$ in its closure $\{ M, G \}$. Since $\mathrm{dim}_\field(U[1], M) = 1$, we know $(U[1], -)$ is given by an $M$-minimal pp-pair by \cite[5.3.10]{prest}.

Dually, each adic is isolated in its closure by a minimal pp-pair.

Finally, as the generic module $G$ has finite endo-length, there is a $G$-minimal pp-pair which, trivially, isolates $G$ in its closure $\{ G \}$ (specifically, choose pp formula $\phi$ such that the length of $\phi(G)$ (as an $\mathrm{End}_R(G)$-module) is minimal but non-zero, then $\phi/(x=0)$ is $G$-minimal). 
\end{proof}
\end{proposition}

\begin{corollary}\label{cor:ASFullSupp}
Let $A$ be a tame hereditary algebra with $S$ a simple regular $A$-module. Given integers $0 \le m \le n$, let $\mathbb{T}$ be the tube of $\modR{A[S,n,m]}$ containing $S$, with $\mathfrak{C} = \mathrm{cl}(\mathbb{T})$ its Ziegler closure, then
\begin{enumerate}[(i)]
\item The generic module in $\mathfrak{C}$ is an $A$-module.
\item Every pr\"ufer module in $\mathfrak{C}$ is an $A[S, n]$-module.
\item Every adic module in $\mathfrak{C}$ is an $A[S,n,m]/I$-module, where $I$ is the ideal generated by the extension vertices $\omega_1, \dots, \omega_n$.
\end{enumerate}
In particular, if $m \ge 1$, then there is no infinite-dimensional module in $\mathrm{cl}(\mathbb{T})$ with full support on the quiver of $A[S,n,m]$.
\begin{proof}
This follows from the construction in the proof of Theorem~\ref{prop:ASTubeClosure}. Each pr\"ufer module in $\mathfrak{C}$ is either an $A$-module or is the image of an $A$-module after a sequence of one-point extensions. Similarly, each adic module in $\mathfrak{C}$ is either an $A$-module or the image of an $A$-module after a sequence of one-point coextensions.
\end{proof}
\end{corollary}

\begin{example}\label{ex:ASKroncker}
Let $C$ be the tame Kronecker algebra and $S$ the simple regular module $\xymatrix{\field & \field \ar@/_/[l]|{\lambda} \ar@/^/[l]|1}$ with $\lambda \in \field$. We know $S$ lies in a homogenous (rank $t=1$) tube of $\modR{C}$. One can calculate that $C[S,2,2]$ is given by the following quiver with the relations $\alpha_2\beta_1 = \lambda\alpha_2\beta_2$, $\beta_1\gamma_1=\lambda\beta_2\gamma_1$, and $\alpha_1\alpha_2\beta_1\gamma_1\gamma_2 = 0$.
\begin{equation*}
\xymatrix{
a & \circ \ar[l]_{\gamma_2} & \circ \ar[l]_{\gamma_1} & \circ \ar@/_/[l]_{\beta_1} \ar@/^/[l]^{\beta_2} & \circ \ar[l]_{\alpha_2} & b \ar[l]_{\alpha_1}
}
\end{equation*}
The module $S$ lies in a tube of the form $\Gamma(0,2,2)$ in $\modR{C[S,2,2]}$. This tube is depicted below with dimension vectors of corresponding modules displayed.
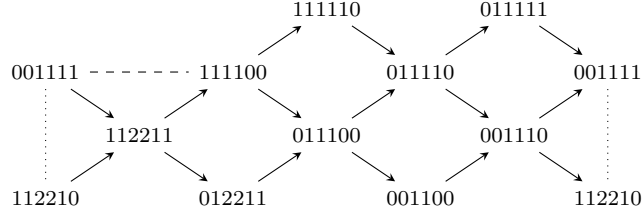
\begin{equation*}
\begin{tikzpicture}[every node/.style={anchor=center, font=\scriptsize},>=stealth]
\matrix [matrix of math nodes, row sep=0.4cm, column sep=0.08cm, anchor=south east] (a) at (0,0) {
&&& 111110 && 011111 \\
001111 && 111100 && 011110 && 001111 \\
& 112211 && 011100 && 001110 \\
112210 && 012211 && 001100 && 112210 \\
};
\draw[->] (a-2-3) to (a-1-4);
\draw[->] (a-1-4) to (a-2-5);
\draw[->] (a-2-5) to (a-1-6);
\draw[->] (a-1-6) to (a-2-7);
\draw[->] (a-2-1) to (a-3-2);
\draw[->] (a-3-2) to (a-2-3);
\draw[->] (a-2-3) to (a-3-4);
\draw[->] (a-3-4) to (a-2-5);
\draw[->] (a-2-5) to (a-3-6);
\draw[->] (a-3-6) to (a-2-7);
\draw[->] (a-4-1) to (a-3-2);
\draw[->] (a-3-2) to (a-4-3);
\draw[->] (a-4-3) to (a-3-4);
\draw[->] (a-3-4) to (a-4-5);
\draw[->] (a-4-5) to (a-3-6);
\draw[->] (a-3-6) to (a-4-7);
\draw[dashed] (a-2-1) to (a-2-3);
\draw[dotted] (a-2-1) to (a-4-1);
\draw[dotted] (a-2-7) to (a-4-7);
\end{tikzpicture}
\end{equation*}
The module $S$ has dimension vector $001100$. The modules with dimension vectors $011111$ and $111110$ are the projective-injective modules corresponding to vertices $a$ and $b$. The closure of this tube in $\Zg_{C[S,2,2]}$ contains the 3 pr\"ufer modules given by direct limits along the rays beginning at $111100$, $111110$, and $011111$. By the placement of $S=001100$, the $111100$-pr\"ufer is the $S$-pr\"ufer of $\Zg_C$. The closure contains the 3 adic modules given by inverse limits along the corays ending at $001111$, $111110$, and $011111$. The $001111$-adic is the $S$-adic of $\Zg_C$.

\end{example}

\begin{example}\label{ex:ASCanonical}
Continuing Example~\ref{ex:CanonicalTubular}, the quiver of $A$ is a Euclidean quiver of type $\mathbb{E}_6$ and $A$ is a tame hereditary algebra containing a tubular family of type $(2,3,3)$. The module $X$ is simple regular and lies in the unique tube of rank 2. The one-point extension $A[X]$ is the canonical tubular algebra of tubular type $(3,3,3)$. 
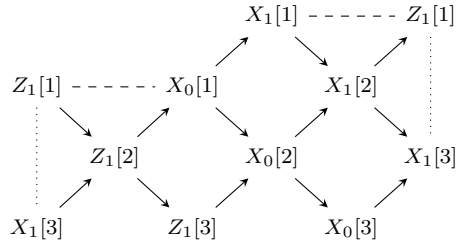
\begin{equation*}
\begin{tikzpicture}[every node/.style={anchor=center, font=\scriptsize},>=stealth]
\matrix [matrix of math nodes, row sep=0.4cm, column sep=0.08cm, anchor=south east] (a) at (0,0) {
&&& X_1[1] && Z_1[1] \\
Z_1[1] && X_0[1] && X_1[2] \\
& Z_1[2] && X_0[2] && X_1[3] \\
X_1[3] && Z_1[3] && X_0[3] \\
};
\draw[->] (a-2-3) to (a-1-4);
\draw[->] (a-1-4) to (a-2-5);
\draw[->] (a-2-5) to (a-1-6);
\draw[->] (a-2-1) to (a-3-2);
\draw[->] (a-3-2) to (a-2-3);
\draw[->] (a-2-3) to (a-3-4);
\draw[->] (a-3-4) to (a-2-5);
\draw[->] (a-2-5) to (a-3-6);
\draw[->] (a-4-1) to (a-3-2);
\draw[->] (a-3-2) to (a-4-3);
\draw[->] (a-4-3) to (a-3-4);
\draw[->] (a-3-4) to (a-4-5);
\draw[->] (a-4-5) to (a-3-6);
\draw[dashed] (a-1-4) to (a-1-6);
\draw[dashed] (a-2-1) to (a-2-3);
\draw[dotted] (a-2-1) to (a-4-1);
\draw[dotted] (a-1-6) to (a-3-6);
\end{tikzpicture}
\end{equation*}
In the notation of this section $A[X] = A[X,1,0]$ and $X$ lies in a tube of $\modR{A[X]}$ having the form $\Gamma(2,1,0)$ drawn above (with $X$ corresponding to $X_0[1]$). The algebra $A[X,1,1]$ is given by the following quiver with relations $\alpha_3\epsilon = 0$, $\beta_3\delta=\beta_3\epsilon$, $\gamma_1\delta=0$, and $\alpha_1\alpha_2\alpha_3 + \beta_1\beta_2\beta_3+\gamma_1\gamma_2\gamma_3 = 0$.
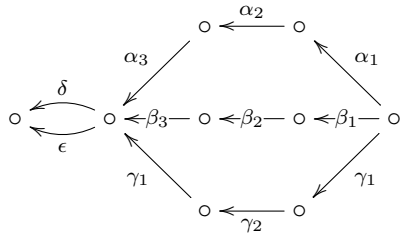
\begin{equation*}
\xymatrix{
&& \circ \ar[dl]_{\alpha_3} & \circ \ar[l]_{\alpha_2} \\
\circ & \circ \ar@/_/[l]_{\delta} \ar@/^/[l]^{\epsilon} & \circ \ar[l]|{\beta_3} & \circ \ar[l]|{\beta_2} & \circ \ar[l]|{\beta_1} \ar[ul]_{\alpha_1} \ar[dl]^{\gamma_1}  \\ 
&& \circ \ar[ul]^{\gamma_1} & \circ \ar[l]^{\gamma_2}
}
\end{equation*}
The module $X$ lies in a tube of $\modR{A[X,1,1]}$ having the form $\Gamma(2,1,1)$ drawn below (beware, as coray insertion doesn't leave rays intact -- in that additional modules will be inserted into rays -- the modules from the previous diagram are renamed --- $X$ now corresponds to $X_0[2]$).
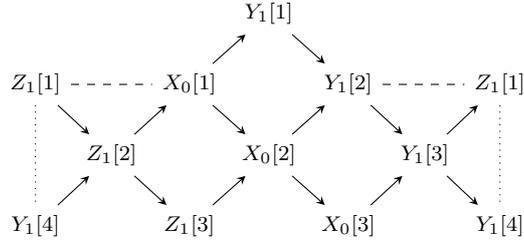
\begin{equation*}
\begin{tikzpicture}[every node/.style={anchor=center, font=\scriptsize},>=stealth]
\matrix [matrix of math nodes, row sep=0.4cm, column sep=0.08cm, anchor=south east] (a) at (0,0) {
&&& Y_1[1] \\
Z_1[1] && X_0[1] && Y_1[2] && Z_1[1] \\
& Z_1[2] && X_0[2] && Y_1[3] \\
Y_1[4] && Z_1[3] && X_0[3] && Y_1[4] \\
};
\draw[->] (a-2-3) to (a-1-4);
\draw[->] (a-1-4) to (a-2-5);
\draw[->] (a-2-1) to (a-3-2);
\draw[->] (a-3-2) to (a-2-3);
\draw[->] (a-2-3) to (a-3-4);
\draw[->] (a-3-4) to (a-2-5);
\draw[->] (a-2-5) to (a-3-6);
\draw[->] (a-3-6) to (a-2-7);
\draw[->] (a-4-1) to (a-3-2);
\draw[->] (a-3-2) to (a-4-3);
\draw[->] (a-4-3) to (a-3-4);
\draw[->] (a-3-4) to (a-4-5);
\draw[->] (a-4-5) to (a-3-6);
\draw[->] (a-3-6) to (a-4-7);
\draw[dashed] (a-2-1) to (a-2-3);
\draw[dashed] (a-2-5) to (a-2-7);
\draw[dotted] (a-2-1) to (a-4-1);
\draw[dotted] (a-2-7) to (a-4-7);
\end{tikzpicture}
\end{equation*}
The closure of this tube in $\Zg_{A[X,1,1]}$ contains: the 3 pr\"ufer modules of the rays beginning $Z_1[1]$, $X_0[1]$, and $Y_1[1]$; the 3 adic modules of the corays ending at $Z_1[1]$, $Y_1[1]$, and $Y_1[2]$; and the unique generic module. The $Y_1$-pr\"ufer and $Y_1$-adic are the only infinite-dimensional points that are not $A$-modules (the pr\"ufer and adic at $Z_1[1]$ are the $\tau X$-pr\"ufer and $\tau X$-adic respectively, and the pr\"ufer at $X_0[1]$ and the adic at $Y_1[2]$ are respectively the $X$-pr\"ufer and $X$-adic). 
\end{example}

\section{Trivial Extensions of Canonical Tubular Algebras}\label{sec:TrivExtCanonical}

The \textbf{canonical algebras} introducted by Ringel \cite[3.7]{ringel} are particular one-point extensions of hereditary algebras. The canonical algebras of ``tubular type'' \cite[5]{ringel}, or \textbf{canonical tubular algebras}, form a class of simply connected algebras whose trivial extensions are non-domestic but of polynomial growth \cite{nehring_skowronski}.

Let $R$ be a canonical algebra of tubular type $\bar{n} = (n_1, \dots, n_r)$ (i.e. one of $(2,2,2,2)$, $(3,3,3)$, $(2,4,4)$, or $(2,3,6)$) and let $T = R \ltimes R^\star$ be its trivial extension. For $\bar{n} = (3,3,3)$ see Examples~\ref{ex:CanonicalTubular} and \ref{ex:CanAlgTrivExt}, which we'll use as a running example throughout. The finite-dimensional representation theory of $T$ is obtained from that of $\hat{R}$ in the following way.

\begin{example}\label{ex:CanAlgTrivExt}
The trivial extension $T = R \ltimes R^\star$ for $\bar{n} = (3,3,3)$.
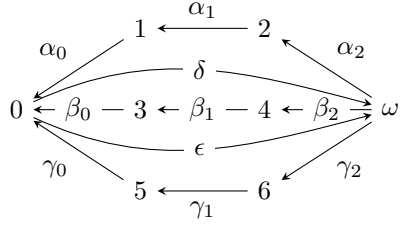
\begin{equation*}
\begin{tikzpicture}[every node/.style={anchor=center, font=\small},>=stealth]
\matrix [matrix of math nodes, row sep=0.6cm, column sep=1.2cm, anchor=south east] (R) at (3,0) {
	& 1 & 2 \\
	0 & 3 & 4 & \omega \\
	& 5 & 6  \\
};
\path (R) +(0,-2.0) node {$\alpha_2\alpha_1\alpha_0 + \beta_2\beta_1\beta_0 + \gamma_2\gamma_1\gamma_0 = 0,$};
\path (R) +(0,-2.5) node {$\delta\gamma_2 = 0,\,\epsilon\alpha_2=0,\,\delta\beta_2=\epsilon\beta_2,$};
\path (R) +(0,-3.0) node {$\alpha_0\epsilon = 0,\,\gamma_0\delta=0,\,\beta_0\delta=\beta_0\epsilon,$};
\path (R) +(0,-3.5) node {$\alpha_1\alpha_0\delta\alpha_2\alpha_1 = 0,\,\beta_1\beta_0\delta\beta_2\beta_1 = 0,\,\gamma_1\gamma_0\epsilon\gamma_2\gamma_1 = 0$.};
\draw[->] (R-1-2) to node[above left] {$\alpha_0$} (R-2-1);
\draw[->] (R-1-3) to node[above] {$\alpha_1$} (R-1-2);
\draw[->] (R-2-4) to node[above right] {$\alpha_2$} (R-1-3);
\draw[->] (R-2-2) to node[midway, fill=white] {$\beta_0$} (R-2-1);
\draw[->] (R-2-3) to node[midway, fill=white] {$\beta_1$} (R-2-2);
\draw[->] (R-2-4) to node[midway, fill=white] {$\beta_2$} (R-2-3);
\draw[->] (R-3-2) to node[below left] {$\gamma_0$} (R-2-1);
\draw[->] (R-3-3) to node[below] {$\gamma_1$} (R-3-2);
\draw[->] (R-2-4) to node[below right] {$\gamma_2$} (R-3-3);
\draw[->] (R-2-1) to[out=25,in=165] node[midway, fill=white] {$\delta$} (R-2-4);
\draw[->] (R-2-1) to[out=-25,in=-165] node[midway, fill=white] {$\epsilon$} (R-2-4);
\end{tikzpicture}
\end{equation*}
\end{example}

Let $\nu: \hat{R} \to \hat{R}$ be the Nakayama automorphism of $\hat{R}$ and $F: \hat{R} \to T$ the Galois covering of Example~\ref{ex:GaloisCoverings}. The induced (restricted) push-down functor $F_\lambda: \modR{\hat{R}} \to \modR{T}$ preserves indecomposable modules and almost-split sequences \cite[3.5, 3.6(a)]{gabriel}. Furthermore, since $\hat{R}$ is ``locally support finite'', the restricted push-down functor is essentially surjective \cite{dowbor_lenzing_skowronski} \cite{happel_ringel} \cite{nehring_skowronski}, and induces a Galois covering $\indR{\hat{R}} \to \indR{T}$ \cite[3.6(c)]{gabriel}. Thus
\begin{equation}\label{eq:indT}
\indR{T} \simeq (\indR{\hat{R}})/\langle\mathrm{\nu}\rangle
\end{equation}
where $\nu$ acts on $\indR{\hat{R}}$ by translation $M \mapsto {}^\nu M$ (i.e.~restriction along $\nu^{-1}$). Also $\Gamma_T \simeq \Gamma_{\hat{R}}/\langle\nu\rangle$ for the induced morphism $\nu: \Gamma_{\hat{R}} \to \Gamma_{\hat{R}}$. In this way, from the structure of $\modR{\hat{R}}$, we obtain the structure of $\modR{T}$. This is done in \cite{happel_ringel} and, for $R$ a general tubular algebra, in \cite{nehring_skowronski}. Further details are contained throughout this section.

\paragraph{}
The quiver $\mcal{Q}_{\hat{R}}$ of $\hat{R}$ is obtained from $\bigsqcup_{i \in \mathbb{Z}} \mcal{Q}_R[i]$ -- a disjoint union of copies of $\mcal{Q}_R$ -- by adding two parallel arrows from the unique sink $x_0[i]$ of $\mcal{Q}_R[i]$ to the unique source $x_\omega[i-1]$ of $\mcal{Q}_R[i-1]$ for each $i \in \mathbb{Z}$. If $\mcal{Q}^\prime$ is any subquiver of $\mcal{Q}_{\hat{R}}$, define $\mcal{Q}^\prime[i] := \nu^i(\mcal{Q}^\prime)$. We will identify $\mcal{Q}_R$ with $\mcal{Q}_R[0]$.

\begin{example}
The quiver $\mcal{Q}_{\hat{R}}$ of $\hat{R}$ for $\bar{n} = (3,3,3)$.
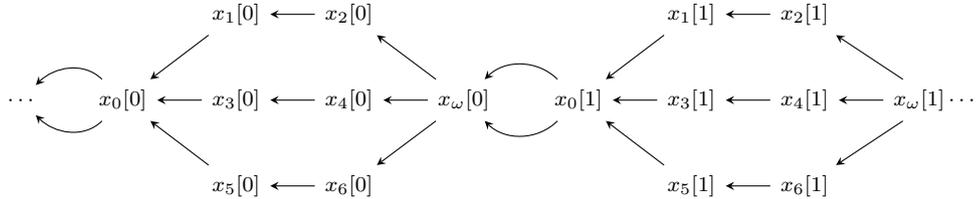
\begin{equation*}
\begin{tikzpicture}[every node/.style={anchor=center, font=\scriptsize},>=stealth]
\matrix [matrix of math nodes, row sep=0.6cm, column sep=0.6cm, anchor=south east] (a) at (3,0) {
	&& x_1[0] & x_2[0] &&& x_1[1] & x_2[1] \\
	\cdots & x_0[0] & x_3[0] & x_4[0] & x_\omega[0] & x_0[1] & x_3[1] & x_4[1] & x_\omega[1] \cdots \\
	&& x_5[0] & x_6[0] &&& x_5[1] & x_6[1]  \\
};
\draw[->] (a-2-2) to[out=135,in=45] (a-2-1);
\draw[->] (a-2-2) to[out=-135,in=-45] (a-2-1);
\draw[->] (a-1-3) to (a-2-2);
\draw[->] (a-1-4) to (a-1-3);
\draw[->] (a-1-7) to (a-2-6);
\draw[->] (a-1-8) to (a-1-7);
\draw[->] (a-2-3) to (a-2-2);
\draw[->] (a-2-4) to (a-2-3);
\draw[->] (a-2-5) to (a-1-4);
\draw[->] (a-2-5) to (a-2-4);
\draw[->] (a-2-5) to (a-3-4);
\draw[->] (a-2-6) to[out=135,in=45] (a-2-5);
\draw[->] (a-2-6) to[out=-135,in=-45] (a-2-5);
\draw[->] (a-2-7) to (a-2-6);
\draw[->] (a-2-8) to (a-2-7);
\draw[->] (a-2-9) to (a-1-8);
\draw[->] (a-2-9) to (a-2-8);
\draw[->] (a-2-9) to (a-3-8);
\draw[->] (a-3-3) to (a-2-2);
\draw[->] (a-3-4) to (a-3-3);
\draw[->] (a-3-7) to (a-2-6);
\draw[->] (a-3-8) to (a-3-7);
\end{tikzpicture}
\end{equation*}
\end{example}

Define the following subquivers of $\mcal{Q}_{\hat{R}}$. If just a set of vertices is given, then we mean the full subquiver determined by those vertices.
\begin{itemize}
\item $\mcal{Q}_{C_0} = \mcal{Q}_R - \{ x_\omega[0] \}$,
\item $\mcal{Q}_{C_1} = \mcal{Q}_R - \{ x_0[0] \}$,
\item $\mcal{Q}_{C_2} = \{ x_\omega[0],\, x_0[1] \}$,
\item $\mcal{Q}_{D_0} = \mcal{Q}_{C_0} \cup \mcal{Q}_{C_1}$,
\item $\mcal{Q}_{D_1} = \mcal{Q}_{C_1} \cup \mcal{Q}_{C_2}$,
\item $\mcal{Q}_{D_2} = \mcal{Q}_{C_2} \cup \mcal{Q}_{C_0}[1]$,
\item $\mcal{Q}_{E_1} = \mcal{Q}_{D_0} \cup \mcal{Q}_{D_1}$,
\item $\mcal{Q}_{E_2} = \mcal{Q}_{D_1} \cup \mcal{Q}_{D_2}$,
\item $\mcal{Q}_{E_0} = \mcal{Q}_{D_2} \cup \mcal{Q}_{D_0}[1]$.
\end{itemize}
Each of the above quivers $\mcal{Q}_A$ (with the induced relations) determines a corresponding finite-dimensional algebra $A$. Moreover, $\mcal{Q}_A$ is path-complete as a subquiver of $\mcal{Q}_{\hat{R}}$, so the category of $A$ can be identified with the full convex subcategory of $\hat{R}$ determined by the vertices of $\mcal{Q}_A$. The push-down functor $F_\lambda: \ModR{\hat{R}} \to \ModR{T}$ then restricts to an interpretation functor $F_\lambda|_{\ModR{A}}: \ModR{A} \to \ModR{T}$ by Corollary~\ref{cor:ConvexRes} which full and faithful in case $A = D_i$.

Define $C_i$ for all $i \in \mathbb{Z}\backslash\{0,1,2\}$ by $C_i := \nu^j(C_k)$ if $i = 3j + k$ for $k \in \{ 0,1,2 \}$ and $j \in \mathbb{Z}$. Make analogous definitions for $D_i$ and $E_i$.

\begin{proposition}\label{prop:HatRSubCat}
The following facts hold:
\begin{enumerate}[(i)]
\item The algebras $C_i$ are tame hereditary and comprise all tame concealed full convex subcategories of $\hat{R}$.
\item The algebras $D_i$ are tubular: being a tubular extension of $C_i$ and a tubular coextension of $C_{i+1}$. They comprise all tubular full subcategories of $\hat{R}$. Moreover, $\hat{R} \simeq \hat{D_i}$ and $T \simeq D_i \ltimes D_i$.
\item Every finite-dimensional $\hat{R}$-module is an $E_i$-module for some $i \in \mathbb{Z}$.
\end{enumerate}
\begin{proof}
See \cite{happel_ringel} and \cite{nehring_skowronski} (also \cite[3.1]{skowronski3}). We point out a few details. $D_0 = R$ is the canonical tubular algebra: a one-point extension of $C_0$ (a tame hereditary of type $\widetilde{\mathbb{D}}_n$ or $\widetilde{\mathbb{E}}_n$ in ``subspace'' orientation) and a one-point coextension of $C_1$ (a tame hereditary of type $\widetilde{\mathbb{D}}_n$ or $\widetilde{\mathbb{E}}_n$ in ``factorspace'' orientation). $C_2$ is the tame Kronecker algebra. $D_1$ and $D_2$ are the so-called ``squid algebras'' appearing in \cite[2]{reiten_ringel}.
\end{proof}
\end{proposition}

For $i \in \mathbb{Z}$ write
\begin{equation}\label{eq:modCi}
\modR{C_i} = \mcal{P}_i \vee \mcal{R}_i \vee \mcal{Q}_i
\end{equation}
where $\mcal{P}_i$ is the preprojective component, $\mcal{Q}_i$ the preinjective component, and $\mcal{R}_i$ the regular component. Similarly, for $i \in \mathbb{Z}$ write
\begin{equation}\label{eq:modDi}
\modR{D_i} = \mcal{P}_i \vee \mcal{T}^0_i \vee \mcal{M}_{i,i+1} \vee \mcal{T}^\infty_{i+1} \vee \mcal{Q}_{i+1}
\end{equation}
where $\mcal{T}^0_i$ is the tubular family of slope 0 (it is obtained from $\mcal{R}_i$ by ray insertion); $\mcal{T}^\infty_{i+1}$ is the tubular family of slope $\infty$ (it is obtained from $\mcal{R}_{i+1}$ by coray insertion); and $\mcal{M}_{i,i+1} = \bigvee_{q \in \mathbb{Q}^+} \mcal{T}^q_i$ are the stable tubular families of positive rational slope. 

\begin{lemma}\label{lemma:modEi}
For $i \in \mathbb{Z}$ we have
\begin{equation}\label{eq:modEi}
\modR{E_i} = \mcal{P}_{i-1} \vee \mcal{T}^0_{i-1} \vee \mcal{M}_{i-1,i} \vee \mcal{T}_i \vee \mcal{M}_{i,i+1} \vee \mcal{T}^\infty_{i+1} \vee \mcal{Q}_{i+1}
\end{equation}
where $\mcal{T}_i$ is a standard non-stable separating tubular $\mathbb{P}^1(\field)$-family containing all regular $C_i$-modules.
\begin{proof}
See \cite[\S 3]{nehring_skowronski}. The structure of the tubular family $\mcal{T}_i$ is described further in Propositions~\ref{prop:ZgE1}---\ref{prop:ZgE2} below. The remaining components are those of $\modR{D_{i-1}}$ and $\modR{D_i}$ from \eqref{eq:modDi}.
\end{proof}
\end{lemma}

\begin{lemma}\label{lemma:modRhat}
We have
\begin{equation}\label{eq:modRhat}
\modR{\hat{R}} \simeq \bigvee_{i \in \mathbb{Z}} \left( \mcal{M}_{i,i+1} \vee \mcal{T}_i \right)
\end{equation}
\begin{proof}
See \cite{happel_ringel} and \cite{nehring_skowronski}.
\end{proof}
\end{lemma}

For $i \in \mathbb{Z}$ we will describe the Ziegler spectrum $\Zg_{E_i}$ in so far as providing a finite cover of closed sets whose relative (subspace) topologies are known. Specifically, we give a family of closed subsets $\mathfrak{C}_1, \dots, \mathfrak{C}_m \subseteq \Zg_{E_i}$ and write
\begin{equation*}
\Zg_{E_i} \simeq \bigcup_{j=1,\dots,m} \mathfrak{C}_j
\end{equation*}
to mean an equality (of sets) such that a subset $\mathfrak{C} \subseteq \Zg_{E_i}$ is closed if and only if $\mathfrak{C} \cap \mathfrak{C}_j$ is closed in $\mathfrak{C}_j$, for each $j = 1, \dots, m$.

We have canonical projections $E_i \to D_j$ which induce closed embeddings $\Zg_{D_j} \to \Zg_{E_i}$, for $j=i-1$ and $j=i$, by Proposition~\ref{prop:InterpFunc}. In this way, we consider $\Zg_{D_{i-1}}$ and $\Zg_{D_i}$ as closed subsets of $\Zg_{E_i}$. The complement $\Zg_{E_i} \backslash (\Zg_{D_{i-1}} \cup \Zg_{D_i})$ then consists precisely of the points corresponding to $E_i$-modules supported over the full quiver $\mcal{Q}_{E_i}$. The intersection $\Zg_{D_{i-1}} \cap \Zg_{D_i}$ may be identified with $\Zg_{C_i}$. In the following results, we extend $\{ \mathfrak{C}_1 = \Zg_{D_{i-1}}$, $\mathfrak{C}_2 = \Zg_{D_i} \}$ to the desired cover of $\Zg_{E_i}$ by adding the closure of finitely many tubes. Each such tube will be one of the non-stable tubes introduced in Section~\ref{sec:NonStableTubes} and whose closure is described in Theorem~\ref{prop:ASTubeClosure} and Proposition~\ref{prop:ASTubeTop}.

\begin{proposition}\label{prop:ZgE1}
There exists a simple regular $C_1$-module $S$ lying in a tube $\mathbb{T} \in \mcal{T}_1$ of $\modR{E_1}$, such that $E_1 = C_1[S,1,1]$ and
\begin{equation*}
\Zg_{E_1} \simeq \Zg_{D_0} \,\cup\, \mathrm{cl}(\mathbb{T}) \,\cup\, \Zg_{D_1}
\end{equation*}
Every infinite-dimensional point in $\Zg_{E_1}$ is either a $D_0$- or a $D_1$-module.
\begin{proof}
The module $S$ is such that $D_0 = [S]C_1$ and $D_1 = C_1[S]$. By Lemma~\ref{lemma:modEi} and Proposition~\ref{prop:AStube}, the tube $\mathbb{T}$ -- containing $S$ -- lies in the tubular family $\mcal{T}_i$ of $\modR{E_i}$, since this is the tubular family containing all regular $C_1$-modules. If $\omega$ (resp. $\sigma$) denotes the extension (resp. coextension) vertex of $E_1 = C_1[S,1,1]$ over $C_1$, then $D_0 = E_1/\langle\omega\rangle$ and $D_1 = E_1/\langle\sigma\rangle$.

It follows that the closed subset $\Zg_{D_0} \cup\, \mathrm{cl}(\mathbb{T}) \cup\, \Zg_{D_1} \subseteq \Zg_{E_1}$ contains all finite-dimensional points. We have equality by the density of these points, by Proposition~\ref{prop:ZgFdAlg}. The final claim follows from Corollary~\ref{cor:ASFullSupp}, since any indecomposable module not supported entirely on $\mcal{Q}_{E_1}$ has support contained within one of the subquivers $\mcal{Q}_{D_0}$ or $\mcal{Q}_{D_1}$.
\end{proof}
\end{proposition}

\begin{proposition}\label{prop:ZgE0}
There exists a simple regular $C_0$-module $S$ lying in a tube $\mathbb{T} \in \mcal{T}_0$ of $\modR{E_0}$, such that $E_0 = C_0[S,1,1]$ and
\begin{equation*}
\Zg_{E_0} \simeq \Zg_{D_0} \,\cup\, \mathrm{cl}(\mathbb{T}) \,\cup\, \Zg_{D_2}
\end{equation*}
Every infinite-dimensional point of $\Zg_{E_0}$ is either a $D_0$- or a $D_2$-module.
\begin{proof}
The proof is similar to Proposition~\ref{prop:ZgE1}. In fact $E_0 \simeq (E_1)^{\op}$ and the argument is essentially dual.
\end{proof}
\end{proposition}

\begin{proposition}\label{prop:ZgE2}
If $\bar{n} = (n_1, n_2, \dots, n_r)$, then there exist simple regular modules $S_1, \dots, S_r \in \modR{C_1}$ and tubes $\mathbb{T}_i \in \mcal{T}_2$ of $\modR{E_2}$ containing $S_i$, for $i=1,\dots,r$ respectively, such that
\begin{equation*}
\Zg_{E_2} \simeq \Zg_{D_1} \,\cup\, \left( \bigcup_{i=1}^r \mathrm{cl}(\mathbb{T}_i) \right) \,\cup\, \Zg_{D_2}
\end{equation*}
Every infinite-dimensional point of $\Zg_{E_2}$ is either a $D_1$- or a $D_2$-module.
\begin{proof}
Note $D_2$ is the tubular extension $C_2[S_1, B_1]\cdots[S_r, B_r]$ of $C_2$ where $B_i$ is a subspace branch of length $n_i - 1$. The modules $S_i$ are simple regular $C_2$-modules lying in distinct homogenous tubes $\mathbb{R}_i \in \mcal{R}_2$ \cite[\S 2, p.8]{reiten_ringel}. The tubular family $\mcal{T}^0_2$ of $\modR{D_2}$ is obtained from $\mcal{R}_2$ (the regular $C_2$-modules) by $n_i-1$ consecutive ray insertions in $\mathbb{R}_i$ -- beginning at $S_i$ -- for each $i=1,\dots,r$. By Proposition~\ref{prop:AStube}, the resulting tube $\mathbb{S}_i$ of $\modR{D_2}$ containing $S_i$ is of the form $\Gamma(0, n_i, 0)$ with $S_i$ corresponding to vertex $X_0[1]$. To obtain $E_2$ from $D_2$, and the tubes $\mathbb{T}_i$ for each $i=1,\dots,r$, we must perform $n_i-1$ consecutive one-point coextensions, corresponding to $n_i-1$ consecutive coray insertions in each $\mathbb{S}_i$. The resulting tubes $\mathbb{T}_i$ lie in the tubular family $\mcal{T}_2$ of $\modR{E_2}$ obtained in this way from $\mcal{T}^0_2$ --- the remaining tubes are left intact and consist of regular $C_2$-modules.

By Lemma~\ref{lemma:modEi}, the closed subset $\Zg_{D_1} \cup \Zg_{D_2} \subseteq \Zg_{E_1}$ contains all finite-dimensional points except for some (in fact, infinitely many) in each of the tubes $\mathbb{T}_1, \dots, \mathbb{T}_r$. Thus the claim follows from the density of the finite-dimensional points by Proposition~\ref{prop:ZgFdAlg}.

In the above construction, we have built each of the tubes $\mathbb{S}_1, \dots, \mathbb{S}_r$ before constructing any of the tubes $\mathbb{T}_1, \dots, \mathbb{T}_r$. We can go about things differently and construct $\mathbb{T}_1$ say, from $\mathbb{S}_1$, before constructing $\mathbb{S}_2, \dots, \mathbb{S}_r$ as follows.

Let $A_0 := C_2[S_1, B_1]$, that is, $A_0 = C_2[S_1, n_1-1]$ --- recall $n_1-1$ is the length of the branch $B_1$. If $\mathbb{A}_0$ is the tube of $\modR{A_0}$ containing $S_1$, then $\mathbb{A}_0 \simeq \mathbb{S}_1$ via the embedding $\modR{A_0} \to \modR{D_2}$ given by restriction along the projection $D_2 \to A_0$. Then the sequence of algebras $A_j = C_2[S_1, n_1-1, j]$ for $j = 1, \dots, n_1-1$, corresponds to $n_1-1$ consecutive coray insertions. If $\mathbb{A}_{n_1-1}$ is the tube of $\modR{A_{n_1-1}}$ containing $S_1$, then $\mathbb{A}_{n_1-1} \simeq \mathbb{T}_1$ and Theorem~\ref{prop:ASTubeClosure} and Proposition~\ref{prop:ASTubeTop} describe the set $\mathrm{cl}(\mathbb{T}_1)$ via the closed embedding $\Zg_{A_{n_1-1}} \to \Zg_{E_2}$ given by restriction along the projection $E_2 \to A_{n_1-1}$.

Of course, we can repeat the above paragraph for each $i=2,\dots,r$. The final claim of the proposition follows from Corollary~\ref{cor:ASFullSupp} (and the proof thereof).
\end{proof}
\end{proposition}

\begin{example}
In the case $\bar{n} = (3,3,3)$, the algebra $E_0$ is given in Example~\ref{ex:ASCanonical}. The quotients of $E_2$, defined in the proof of Proposition~\ref{prop:ZgE2}, are (up to isomorphism) given in Example~\ref{ex:ASKroncker}.
\end{example}

Recall from Proposition~\ref{prop:HatRSubCat} that $T \simeq D_i \ltimes D^\star_i$ for all $i \in \mathbb{Z}$.


For $i \in \mathbb{Z}$, let $\mcal{D}_i$ denote the definable subcategory of $\ModR{D_i}$ generated by the collection of tubular families $\mcal{M}_{i,i+1}$ in $\modR{D_i}$. Let $\mcal{E}_i$ denote the definable subcategory of $\ModR{E_i}$ generated by the tubular family $\mcal{T}_i$ in $\modR{E}_i$. Under the $\nu$-translation on $\ModR{\hat{R}}$, we have ${}^\nu \mcal{D}_i = \mcal{D}_{i+3}$ and ${}^\nu \mcal{E}_i = \mcal{E}_{i+3}$.

\begin{theorem}\label{theorem:ZgTrivCanonical}
Let $R$ be a canonical tubular algebra and $F: \hat{R} \to T$ the Galois covering of the trivial extension $T = R \ltimes R^\star$. The push-down functor $F_\lambda: \ModR{\hat{R}} \to \ModR{T}$ induces closed embeddings $\Zg(\mcal{D}_i) \to \Zg_T$ and $\Zg(\mcal{E}_i) \to \Zg_T$ such that
\begin{equation}\label{eq:ZgT}
\Zg_T \simeq \bigcup_{i=0,1,2} \Zg(\mcal{D}_i) \,\cup\, \Zg(\mcal{E}_i)
\end{equation}
for the definable subcategories $\mcal{D}_i$ and $\mcal{E}_i$ of $\ModR{\hat{R}}$ defined above. Every finite-dimensional point of $\Zg_T$ is of the form $F_\lambda M$ for some $M \in \mcal{E}_i$. Every infinite-dimensional point of $\Zg_T$ is of the form $F_\lambda M$ for some $M \in \mcal{D}_i$.
\begin{proof}
As noted above, each of $D_i$ and $E_i$ are finite convex full subcategories of $\hat{R}$ (as their respective quivers are path-complete subquivers of $\mcal{Q}_{\hat{R}}$). Hence, for all $i \in \mathbb{Z}$, the restriction of $F_\lambda$ to $\ModR{D_i}$ or $\ModR{E_i}$ is an interpretation functor by Corollary~\ref{cor:ConvexRes} and moreover $F_\lambda|_{\ModR{D_i}}$ is full and faithful. As $F_\lambda$ preserves indecomposability of all (pure-injective) $D_i$-modules, it induces a closed embedding $\Zg_{D_i} \to \Zg_T$ (see Proposition~\ref{prop:InterpFunc}) which restricts to $\Zg(\mcal{D}_i) \to \Zg_T$.

By Propositions~\ref{prop:ZgE1}---\ref{prop:ZgE2} every infinite-dimensional point of $\Zg_{E_i}$ is actually a $D_j$-module for some $j$, so indecomposability is preserved by $F_\lambda$ for such points. As previous noted, $F_\lambda$ preserves indecomposability of all finite-dimensional points \cite[3.5]{gabriel}.  Hence, $F_\lambda|_{\ModR{E_i}}: \ModR{E_i} \to \ModR{T}$ induces a closed and continuous map $\Zg_{E_i} \to \Zg_T$ whose restriction $\Zg(\mcal{E}_i) \to \Zg_T$ is one-to-one.

By Proposition~\ref{prop:HatRSubCat} and the density of $F_\lambda: \modR{\hat{R}} \to \modR{T}$ -- see equation \eqref{eq:indT} -- the following closed subset contains all finite-dimensional points
\begin{equation*}
\bigcup_{i = 0,1,2} \Zg(\mcal{D}_i) \cup \Zg(\mcal{E}_i) \subseteq \Zg_T
\end{equation*}
with equality following from Proposition~\ref{prop:ZgFdAlg}.
\end{proof}
\end{theorem}

\begin{corollary}\label{cor:TubeClosure}
Let $T$ be the trivial extension of a canonical tubular algebra $R$ and let $\mathbb{T}$ be a tube in $\modR{T}$, then the Ziegler closure $\mathrm{cl}(\mathbb{T})$ of $\mathbb{T}$ contains:
\begin{enumerate}[(i)]
\item The pr\"ufer modules obtained by direct limits along the rays of $\mathbb{T}$.
\item The adic modules obtained by inverse limits along the corays of $\mathbb{T}$.
\item A unique generic module.
\end{enumerate}
Each pr\"ufer and adic module has CB rank 1, and the generic module has CB rank 2, in the closure of the tube. Furthermore, $\mathrm{cl}(\mathbb{T})$ satisfies the isolation condition.
\begin{proof}
For each $\mathbb{T}$ of $\modR{T}$ there is a tube $\mathbb{S}$ of $\modR{\widehat{R}}$ such that $\mathbb{T} = F_\lambda \mathbb{S}$. Furthermore, $F_\lambda$ induces a homeomorphism $\mathrm{cl}(\mathbb{S}) \simeq \mathrm{cl}(\mathbb{T})$ which preserves pr\"ufer and adic modules in the following sense.
\begin{enumerate}[(i)]
\item If $U[1] \to U[2] \to \cdots$ is a maximal ray of $\mathbb{S}$, then
\begin{equation*}
F_\lambda(\varinjlim U[i]) = \varinjlim F_\lambda U[i]
\end{equation*}
is the $F_\lambda U$-pr\"ufer of $\mathbb{T}$.
\item If $\cdots [2]V \to [1]V$ is a maximal coray of $\mathbb{S}$, then
\begin{equation*}
F_\lambda(\varprojlim [i]V) = \varprojlim F_\lambda [i]V
\end{equation*}
is the $F_\lambda V$-adic of $\mathbb{T}$ (we know $F_\lambda$ preserves the inverse limit by Corollary~\ref{cor:PushDownLimits}).
\end{enumerate}
Note every pr\"ufer and adic module of $\mathrm{cl}(\mathbb{T})$ is of this form.

The stated CB ranks follow from Proposition~\ref{prop:ASTubeTop} for non-stable tubes and from \cite[Lemma~51]{harland} for stable tubes. We prove the isolation condition holds for $\mathrm{cl}(\mathbb{T})$ in $\Zg_T$. Just as in $\mathrm{cl}(\mathbb{S})$, the finite-dimensional points in $\mathrm{cl}(\mathbb{T})$ are isolated by minimal pp-pairs. Let $i \in \mathbb{Z}$ be such that $\mathbb{T}$ consists of $E_i$-modules (such $i$ exists by Lemmas~\ref{lemma:modEi}---\ref{lemma:modRhat}).

If $G \in \mathrm{cl}(\mathbb{S})$ is the generic module, then $F_\lambda G$ is the generic module of $\mathrm{cl}(\mathbb{T})$. As $G$ is a $C_i$-module, for $M \in \mathrm{cl}(\mathbb{S}) \backslash \{ G \}$ we have $(M, {}^{\nu^n} G) = 0$ for all $n \ne 0$ (for $M$ is an $E_i$-module and $\mcal{Q}_{E_i}$ intersects $\mcal{Q}_{C_{i+3k}}$ only when $k=0$). Therefore $(F_\lambda M, F_\lambda G) \simeq (M, G) = 0$ by Lemma~\ref{lemma:PushDownBij}. Similarly $(F_\lambda G, F_\lambda M) \simeq (G, M) = 0$.

For $M = U[\infty]$ it follows that $(F_\lambda U[1], -)$ isolates $F_\lambda U[\infty]$ in its closure $\{ F_\lambda U[\infty], F_\lambda G \}$ and is given by a minimal pp-pair (just like $(U[1], -)$ is, in Proposition~\ref{prop:ASTubeTop} / \cite[Lemma~51]{harland}). Similarly, every adic module in $\mathrm{cl}(\mathbb{T})$ is isolated in its closure by a minimal pp-pair. Finally, $F_\lambda G$ is generic and isolated by a minimal pp-pair in its closure in the same manner as $G$ (see proof of Proposition~\ref{prop:ASTubeTop}).
\end{proof}
\end{corollary}

The remaining infinite-dimensional points of $\Zg_T$, not described in the preceding result, are the $D_i$-modules of irrational slope.

\section{Self-Injective Algebras of Tubular Type}\label{sec:ZgSelfInj}
The results of the previous sections, and the knowledge of what forms a self-injective algebra can take (outlined in \cite{skowronski}), are applied in this section to deduce the Ziegler spectrum for a large class of self-injective algebras of polynomial growth.

\paragraph{}
If $R$ is a self-injective algebra, let $\mcal{O} \subseteq \Zg_R$ be the subset corresponding to the projective $R$-modules. Let $\Zg^s_R := \Zg_R \backslash \mcal{O}$, then $\Zg_R = \mcal{O} \,\sqcup\, \Zg^s_R$ as a disjoint union. Let $\smodR{R}$ (resp.~$\sModR{R}$) denote the category obtained as the quotient of $\modR{R}$ (resp.~$\ModR{R}$) by the ideal generated by all morphisms that factor through a projective module. We now relate $\Zg^s_R$ to these stable module categories.

We briefly introduce certain triangulated categories (see for instance \cite[\S 1]{happel} for a definition of \textbf{triangulated category}). A notion of purity for ``compactly generated'' triangulated categories is defined algebraically in \cite{beligiannis} and \cite{krause3}, and in model theoretic terms in \cite{garkusha_prest}. If $\mcal{A}$ is such a category, then its Ziegler spectrum $\Zg(\mcal{A})$ is defined as the set of (isomorphism classes of) indecomposable pure-injective objects. Let $\mcal{C} \subseteq A$ be the full subcategory of compact objects, then $\ModR{\mcal{C}} := (\mcal{C}^{\op}, \Ab)$ is a locally coherent category. The functor $\mcal{A} \to \ModR{\mcal{C}}$ -- defined by $X \mapsto (-, X)|_{\mcal{C}}$ for $X \in \mcal{A}$ -- then induces a bijection $\Zg(\mcal{A}) \simeq \Sp(\ModR{\mcal{C}})$ \cite[1.9]{krause3}. Here $\Sp(-)$ denotes the injective spectrum, which is defined (with Ziegler topology) for locally coherent categories in \cite{krause} and \cite{herzog}.

In particular, for a self-injective algebra $R$, the stable module categories $\smodR{R}$ and $\sModR{R}$ are triangulated \cite[\S I.3]{happel}. As such, $\sModR{R}$ is compactly generated and its subcategory of compact objects may be identified with $\smodR{R}$ \cite[\S 1.5]{krause3} \cite[\S 4.1]{herzog}. 

\begin{proposition}[{\cite[6.1]{garkusha_prest}}]\label{prop:garkusha_prest}
If $R$ is a self-injective algebra, then there is a homeomorphism $\Zg^s_R \simeq \Zg(\sModR{R})$.
\begin{proof}
This is a special case of \cite[6.1]{garkusha_prest}. Note, by \cite[1.16]{krause3}, an (indecomposable) $R$-module is pure-injective in $\ModR{R}$ if and only if it is pure-injective as an object of $\sModR{R}$. The stated homeomorphism is given by this identification of points.
\end{proof}
\end{proposition}

We say two algebras $A$ and $B$ are \textbf{stably equivalent} if there exists an equivalence $\smodR{A} \simeq \smodR{B}$.


\begin{proposition}\label{prop:ZgStabEquiv}
If $F: \smodR{A} \to \smodR{B}$ is a stable equivalence for self-injective algebras $A$ and $B$, then $F$ induces a homeomorphism $H: \Zg^s_A \to \Zg^s_B$ with $H(M) = F(M)$ and $H(\tau_A M) = \tau_B H(M)$ for all non-projective $M \in \indR{A}$.
\begin{proof}
As restriction along $F$ is an equivalence, its left adjoint $F_L$ is also an equivalence \cite[IV.4.1]{maclane}. We have
\begin{equation*}
\Zg(\sModR{A}) = \Sp((\smodR{A})^{\op}, \Ab) \xrightarrow{F_L} \Sp((\smodR{B})^{\op}, \Ab) = \Zg(\sModR{B})
\end{equation*}
which gives a homeomorphism $H: \Zg^s_A \to \Zg^s_B$ by Proposition~\ref{prop:garkusha_prest}. Now $F_L$ extends $F$ in the sense that the following diagram commutes (where the vertical arrows are the Yoneda embeddings).
\begin{equation*}
\xymatrix{
\smodR{A} \ar[r]^{F} \ar[d] & \smodR{B} \ar[d] \\
((\smodR{A})^{\op}, \Ab) \ar[r]^{F_L} & ((\smodR{B})^{\op}, \Ab)
}
\end{equation*}
Indeed, let $\mcal{A} = \smodR{A}$ and $\mcal{B} = \smodR{B}$, then given $M \in \mcal{A}$ we have
\begin{equation*}
F_L (\mcal{A}(-, M))= \mcal{B}(-, F-) \otimes_{\mcal{A}} \mcal{A}(-,M) \simeq \mcal{B}(-, FM)
\end{equation*}
by Proposition~\ref{prop:KanExt}. After identifying $M \in \mcal{A}$ with the corresponding representable functor $\mcal{A}(-, M)$ (and similarly for $FM \in \mcal{B}$), we have $H(M) = F(M)$ as claimed.

Finally, $F$ induces an isomorphism of (valued) quivers $\Gamma^s_A \to \Gamma^s_B$ by \cite[X.1.3]{auslander_smalo}. This is automatically an isomorphism of stable translation quivers in this case, in particular $F$ (hence $H$) commutes with the AR translations.
\end{proof}
\end{proposition}

\subsection{Trivial extensions of tubular algebras}

Two algebras $A$ and $B$ are \textbf{tilting-cotilting equivalent}\footnote{This is an equivalence relation by \cite[4.1.2]{ringel}.} if there exists a series of algebras $A_1, \dots, A_n$, with $A_1 = A$, $A_n = B$, and a tilting \emph{or} cotilting $A_i$-module $M_i$ such that $A_{i+1} = \End_{A_i}(M_i)$, for $i = 1, \dots, n-1$.

\begin{theorem}[{\cite{tachikawa_wakamatsu}}]\label{theorem:TiltCotiltStbEqv}
If $A$ and $B$ are tilting-cotilting equivalent algebras, then $A \ltimes A^\star$ and $B \ltimes B^\star$ are stably equivalent.
\begin{proof}
This is (repeated use of) the main theorem of \cite{tachikawa_wakamatsu}.
\end{proof}
\end{theorem}

\begin{lemma}[{\cite[1]{happel_ringel}}]
Every tubular algebra is tilting-cotilting equivalent to a canonical tubular algebra.
\begin{proof}
The proof of \cite[1]{happel_ringel} -- that every tubular algebra is derived equivalent to a canonical tubular algebra -- contains the proof of this statement.
\end{proof}
\end{lemma}

\begin{corollary}\label{cor:ZgTrivTubular}
If $A$ is a tubular algebra, then there exists a canonical tubular algebra $B$ such that $\Zg_{A \ltimes A^\star} \simeq \Zg_{B \ltimes B^\star}$.
\begin{proof}
Choose a canonical tubular algebra $B$ tilting-cotilting equivalent to $A$. By Theorem~\ref{theorem:TiltCotiltStbEqv}, the trivial extensions $A \ltimes A^\star$ and $B \ltimes B^\star$ are stably equivalent. Under the tilting-cotilting equivalence of $A$ and $B$, the number of indecomposable projective modules is preserved, so the homeomorphism $\Zg^s_{A \ltimes A^\star} \simeq \Zg^s_{B \ltimes B^\star}$ given by Proposition~\ref{prop:ZgStabEquiv} extends to a homeomorphism of the whole spectrum.
\end{proof}
\end{corollary}

\subsection{Standard self-injective algebras}
%

A self-injective algebra $R$ is said to be \textbf{standard} if it admits a Galois covering $F: A \to R$ where $A$ is simply connected. A self-injective algebra of \textbf{Euclidean type} (resp. \textbf{tubular type}) is an algebra of the form $\widehat{A}/G$ where $A$ is a tilted algebra (resp. tubular algebra) and $G$ is an admissible group of automorphisms of $\hat{A}$. Both classes of algebras are standard and of polynomial growth \cite[4.1, 5.2]{skowronski}. They give all such algebras (up to isomorphism) by the following theorem.

\begin{theorem}[{\cite[6.1]{skowronski}}]
If $R$ is a self-injective algebra of (representation-infinite) polynomial growth, then $R$ is standard if and only if $R$ is isomorphic to a self-injective algebra of Euclidean type or a self-injective algebra of tubular type.
\end{theorem}


 An automorphism $\phi: \hat{R} \to \hat{R}$ is said to be \textbf{rigid} if, for all $(i, m) \in I \times \mathbb{Z}$, we have $\phi(i, m) = (j, m)$ for some $j \in I$. Similarly, say $\phi$ is \textbf{positive} (resp. \textbf{strictly positive}) if, for all $(i, m) \in I \times \mathbb{Z}$, we have $\phi(i, m) = (j, n)$ for some $j \in I$ and $n \ge m$ (resp. $n > m$) in $\mathbb{Z}$. For example, the Nakayama automorphism $\nu: \hat{R} \to \hat{R}$ is strictly positive.

For $R$ a tubular algebra, the admissible groups of automorphisms of $\hat{R}$ are restricted to the following form.

\begin{proposition}[{\cite[3.8, 3.9]{skowronski3}}]
If $R$ is a tubular algebra, then there exists a positive automorphism $\phi: \hat{R} \to \hat{R}$ such that
\begin{enumerate}[(i)]
\item there exists a rigid automorphism $\delta$ of $\hat{R}$ and integer $t \ge 1$ such that $\delta\phi^t = \nu$, and
\item if $G$ is an admissible group of automorphisms of $\hat{R}$, then $G$ is an infinite cyclic group generated by $\sigma\phi^s$ for some rigid automorphism $\sigma$ and integer $s \ge 1$.
\end{enumerate}
\end{proposition}

Following \cite[\S 3]{skowronski3} a tubular algebra $R$ is called \textbf{exceptional} if $t > 1$ in the above proposition, otherwise $R$ is \textbf{normal} in which case $\phi = \nu$ (and $\delta = 1$).

A canonical tubular algebra $R$ is normal by \cite[3.2]{skowronski3} and Theorem~\ref{theorem:ZgTrivCanonical} immediately generalises. Let $\mcal{D}_i$ and $\mcal{E}_i$ denote the definable subcategories of $\ModR{\hat{R}}$ as defined in Section~\ref{sec:TrivExtCanonical}.

\begin{corollary}\label{cor:CanCover}
If $R$ is a canonical tubular algebra and $F: \hat{R} \to \hat{R}/G$ a Galois covering functor, then there exists an integer $s \ge 1$ such that
\begin{equation}\label{eq:ZgS}
\Zg_{\hat{R}/G} \simeq \bigcup_{i=0,1,\ldots,3s-1} \Zg(\mcal{D}_i) \cup \Zg(\mcal{E}_i)
\end{equation}
\begin{proof}
Since $R$ is normal, we know $G = \langle\varphi\rangle$ with $\varphi = \sigma\nu^s$ for some rigid automorphism $\sigma$ and integer $s \ge 1$.
The (restricted) push-down functor $F_\lambda: \modR{\hat{R}} \to \modR{(\hat{R}/G)}$ induces a Galois covering $\indR{\hat{R}} \to \indR{\hat{R}/G}$ with ${}^\varphi \mcal{D}_i = \mcal{D}_{i+3s}$ and ${}^\varphi \mcal{E}_i = \mcal{E}_{i+3s}$.
\begin{equation*}
\indR{\hat{R}/G} \simeq (\indR{\hat{R}})/\langle\varphi\rangle
\end{equation*}
Then, as in Theorem~\ref{theorem:ZgTrivCanonical}, using Corollary~\ref{cor:ConvexRes}, the push-down functor $F_\lambda: \ModR{\hat{R}} \to \ModR{\hat{R}/G}$ preserves indecomposability for pure-injective $E_i$-modules and induces closed embeddings $\Zg(\mcal{D}_i) \to \Zg_{\hat{R}/G}$ and $\Zg(\mcal{E}_i) \to \Zg_{\hat{R}/G}$. In this way, the right-hand side of \eqref{eq:ZgS} is a closed subset containing all finite-dimensional points, with equality by Proposition~\ref{prop:ZgFdAlg}.
\end{proof}
\end{corollary}

\subsection{Non-standard self-injective algebras}
If $A$ is a self-injective algebra, then the socle of $A_A$ is a two-sided ideal of $A$ (and equals the socle of ${}_A A$); define $A_S = A/\mathrm{soc}(A)$. Two self-injective algebras $A$ and $B$ are \textbf{socle equivalent} if there exists an isomorphism $A_S \simeq B_S$.

\begin{theorem}[{\cite[6.2]{skowronski}}]
Any non-standard self-injective algebra of polynomial growth is socle equivalent to a unique (up to isomorphism) standard self-injective algebra of polynomial growth.
\end{theorem}

\begin{lemma}[{\cite[8.69]{jensen_lenzing}}]
If $A$ is a self-injective algebra and $M \in \ModR{A}$, then $M = N \oplus E$ where $E$ is injective and $N\mathrm{soc}(A) = 0$.
\end{lemma}

\begin{proposition}\label{prop:SocEquiv}
If $A$ and $B$ are socle-equivalent self-injective algebras, then there exists a homeomorphism $\Zg_A \backslash \mcal{O} \simeq \Zg_B \backslash \mcal{O}^\prime$, where $\mcal{O}$ and $\mcal{O}^\prime$ are the subsets of projective-injective $A$- and $B$- modules, respectively.
\begin{proof}
The canonical projection $\pi: A \to A_S$ induces a closed embedding $\Zg_{A_S} \to \Zg_A$. By the previous lemma, the image of this embedding contains all non-injective finite-dimensional points. Hence, as $A$ is self-injective, $\Zg_A \simeq \Zg_{A_S} \,\sqcup\, \mcal{O}$ by Proposition~\ref{prop:ZgFdAlg}. Likewise $\Zg_B \simeq \Zg_{B_S} \,\sqcup\, \mcal{O}^\prime$. The assumed isomorphism $A_S \simeq B_S$ implies a homeomorphism $\Zg_{A_S} \simeq \Zg_{B_S}$ and the result follows.
\end{proof}
\end{proposition}

\begin{corollary}
If $A$ is a self-injective algebra of polynomial growth, then there exists a standard self-injective algebra $B$ of polynomial growth and a homeomorphism $\Zg_A \backslash \mcal{O} \simeq \Zg_B \backslash \mcal{O}^\prime$, where $\mcal{O}$ and $\mcal{O}^\prime$ are the subsets of projective-injective $A$- and $B$- modules, respectively.
\end{corollary}

\section{Appendix}\label{sec:Appendix}
\subsection{Tensor products and adjoints to restriction}\label{sec:TensorProducts}
If $A$ is a small $\field$-linear category, then there exists bifunctors
\begin{align*}
[-, -]_A: (A, \mcal{V})^{\op} \otimes_\field (A, \mcal{V}) \to \mcal{V} \\
- \otimes_A -: (A^{\op}, \mcal{V}) \otimes_\field (A, \mcal{V}) \to \mcal{V}
\end{align*}
where $\mcal{V} = \ModR{\field}$ and $(A, \mcal{V})$ denotes all $\field$-linear functors $A \to \mcal{V}$. We define $\RMod{A} := (A, \mcal{V})$ and $\ModR{A} := (A^{\op}, \mcal{V})$.

The functor $[-, -]_A$ is just the \textbf{hom-functor} of $\RMod{A}$. That is, given $M, N \in \RMod{A}$, then $[M, N]_A = \Hom_{(A, \mcal{V})}(M, N)$ (and this set is usually abbreviated as $\Hom_A(M, N)$). Note $[-,-]_{A^{\op}}$ is the hom-functor of $\ModR{A}$ but we will typically forgo writing the ``$\op$'' when context is clear.

The functor $- \otimes_A -$ is the \textbf{tensor product} and generalises the tensor product over a $\field$-algebra. This bifunctor is defined in \cite[1]{fisher_newell} and \cite[6]{mitchell} (see also \cite{krause4} and at a higher level of generality \cite[3.1]{kelly}). Dual to the hom-functor, it preserves colimits in both variables, and satisfies the ``co-Yoneda'' isomorphisms $M \otimes_A A(a, -) \simeq M(a)$ and $A(-, a) \otimes_A N \simeq N(a)$ for all $M \in \ModR{A}$ and $N \in \RMod{A}$. Additionally, unlike the hom-functor, there exists the symmetry $M \otimes_A N \simeq N \otimes_{A^{\op}} M$, so we can always avoid writing an ``$\op$'' if desired.

An \textbf{$A$-$B$-bimodule} is a functor $X: A \otimes_\field B^{\op} \to \mcal{V}$. Given such a bimodule there are the various induced functors: $[X, -]_B:\ModR{B} \to \ModR{A}$ and $[-, X]_A: (\RMod{A})^{\op} \to \ModR{B}$; similarly $X \otimes_B -: \RMod{B} \to \RMod{A}$ and $-\otimes_A X:\;\ModR{A} \to \ModR{B}$; and the usual adjunction
\begin{equation*}
[M \otimes_A X, N]_B \simeq [M, [X, N]_B]_A
\end{equation*}
for $M \in \ModR{A}$ and $N \in \ModR{B}$, i.e $- \otimes_A X$ is left adjoint to $[X, -]_B$.

In particular, we always have the $A$-$A$-bimodule $A(-, -): A^{\op} \otimes_\field A \to \mcal{V}$ given by the hom-functor of $A$. If $F: A \to B$ is a $\field$-linear functor, we have the $B$-$A$-bimodule $B(F-, -): A^{\op} \otimes_\field B \to \mcal{V}$; and similarly the $A$-$B$-bimodule $B(-, F-): B^{\op} \otimes_\field A \to \mcal{V}$. These bimodules induce the following functors
\begin{align*}
F_L := B(F-, -) \otimes_A \,?\,: \RMod{A} \to \RMod{B} \\
F_R := [B(-, F-), \,?\,]_A: \RMod{A} \to \RMod{B}
\end{align*}
where $?$ represents the variable in $\RMod{A}$. Let $\mathrm{res}_F: \RMod{B} \to \RMod{A}$ denote restriction along $F$, i.e. the functor defined by $M \mapsto M \circ F$ for $M \in \RMod{A}$. These three functors are related and have the following properties.

\begin{proposition}\label{prop:KanExt}
Suppose $F: A \to B$ is a $\field$-linear functor between small $\field$-linear categories $A$ and $B$.
\begin{enumerate}[(i)]
\item The functor $F_L := B(F-, -) \otimes_A \,?$ is left adjoint to $\mathrm{res}_F$.
\item The functor $F_R := [B(-, F-), \,?\,]_A$ is right adjoint to $\mathrm{res}_F$.
\item The restriction $\mathrm{res}_F$ preserves all limits and colimits.
\item If $F$ is full and surjective on objects, then $\mathrm{res}_F$ is full and faithful.
\item If $A \subseteq B$ is a full subcategory and $F: A \to B$ the inclusion, then both adjoints $F_L$ and $F_R$ are full and faithful.
\item The following are equivalent:
\begin{enumerate}[(a)]
\item $F_L$ commutes with products,
\item $F_L$ is an interpretation functor,
\item $B(F-, b)$ is finitely-presented for all $b \in B$,
\item $\mathrm{res}_F$ restricts to a functor $\modR{B} \to \modR{A}$.
\end{enumerate}
\item The following are equivalent
\begin{enumerate}[(a)]
\item $F_R$ commutes with direct limits,
\item $F_R$ is an interpretation functor,
\item $B(b, F-)$ is finitely-presented for all $b \in B$,
\item $\mathrm{res}_F$ restricts to a functor $\Rmod{B} \to \Rmod{A}$.
\end{enumerate}
\end{enumerate}
\begin{proof}
Parts (i) and (ii) are by the hom-tensor adjunctions, since
\begin{equation*}
[B(F-, -), \,?\,]_B \simeq \mathrm{res}_F(\,?\,) \simeq B(-, F-) \otimes_B \,?\,
\end{equation*}
by the (co-)Yoneda isomorphisms \cite[\S 3 Eq.~(3.10)]{kelly}, cf. \cite[4.5, 4.6]{fisher_newell} and \cite[\S6]{mitchell}. Part (iii) follows from (i) and (ii) since $\mathrm{res}_F$ is itself both a left and a right adjoint. Part (iv) is an easy exercise. Part (v) is by \cite[3.4(e), 3.5(d)]{auslander}.

For (vi) (a) $\Leftrightarrow$ (b) note that, as a left adjoint, $F_L$ commutes with direct limits, and is therefore an interpretation functor if and only if it commutes with products.

For (vi) (a) $\Leftrightarrow$ (c) the functor $F_L = B(F-, -) \otimes_A \,?: \RMod{A} \to \RMod{B}$ commutes with products if and only if, given a set of $A$-modules $\{ M_i \}_{i \in I}$, the canonical morphism
\begin{equation*}
B(F-, -) \otimes_A \left(\prod_{i \in I} M_i\right) \to \prod_{i \in I} \left(B(F-, -) \otimes_A M_i \right)
\end{equation*}
is an isomorphism of $B$-modules \cite[V.4, Ex.5]{maclane}. That is, for all $b \in B$, the components
\begin{equation*}
B(F-, b) \otimes_A \left(\prod_{i \in I} M_i\right) \to \prod_{i \in I} \left(B(F-, b) \otimes_A M_i \right)
\end{equation*}
are $\field$-linear isomorphisms. This condition on $b$ is precisely that $B(F-, b)$ is a finitely-presented (right) $A$-module by \cite[I.13.2]{stenstrom}.

For (vi)(c) $\Leftrightarrow$ (d) note that $B(F-, b) = \mathrm{res}_F(B(-, b))$. If $M \in \ModR{B}$ is finitely-presented, then there exists an exact sequence
\begin{equation*}
\bigoplus_{i=1}^n B(-, x_i) \to \bigoplus_{j=1}^m B(-, y_j) \to M \to 0
\end{equation*}
for some $x_i, y_j \in B$. Applying the right exact $\mathrm{res}_F$ gives an exact sequence
\begin{equation*}
\bigoplus_{i = 1}^n B(F-, x_i) \to \bigoplus_{j=1}^m B(F-, y_j) \to \mathrm{res}_F M \to 0
\end{equation*}
and $\mathrm{res}_F M$ is finitely-presented if each $B(F-, x_i)$ and $B(F-, y_j)$ is finitely-presented, since $\modR{A}$ is closed under cokernels \cite[E.1.16]{prest}. Conversely, if $\mathrm{res}_F$ restricts to a functor $\modR{B} \to \modR{A}$, then as $B(-, b)$ is finitely-presented, so too is $\mathrm{res}_F(B-, b) = B(F-, b)$.

The proof of (vii) is similar. As a right adjoint, $F_R$ commutes with products, so is an interpretation functor if and only if it commutes with direct limits. Now, the functor $F_R := [B(-, F-), \,?\,]_A: \RMod{A} \to \RMod{B}$ commutes with direct limits if and only if, given a directed set of $A$-modules $\{M_i \}_{i \in I}$ and $b \in B$, the canonical morphism
\begin{equation*}
\varinjlim [B(b, F-), M_i] \to [B(b, F-), \varinjlim M_i]_A
\end{equation*}
is a $\field$-linear isomorphism. This condition on $b$ is precisely that $B(b, F-)$ is a finitely-presented (left) $A$-module by \cite[V.3.4]{stenstrom}. The proof of (c) $\Leftrightarrow$ (d) is dual to (iv).
\end{proof}
\end{proposition}

\begin{corollary}\label{cor:AdjunctLift}
Using the notation of Proposition~\ref{prop:KanExt}, if $F$ has a right adjoint $G: B \to A$, then $F_L \simeq \mathrm{res}_G$, and the adjunction $F \dashv G$ lifts to an adjunction $\mathrm{res}_G \dashv \mathrm{res}_F$.
\begin{proof}
By assumption we have an isomorphism $B(F-,-) \simeq A(-, G-)$ of $A$-$B$-bimodules, hence $F_L = B(F-, -) \,\otimes_A \,?\, \simeq A(-,G-) \,\otimes_A \,?\, = \mathrm{res}_G$ as claimed.
\end{proof}
\end{corollary}

\begin{corollary}\label{cor:IteratedRes}
Using the notation of Proposition~\ref{prop:KanExt}, if $F: A \to B$ and $G: B \to C$ are $\field$-linear functors between small $\field$-linear categories, then $(GF)_L = G_L \circ F_L$ and $(GF)_R = G_R \circ F_R$ are the left and right adjoints, respectively, to $\mathrm{res}_{GF}$.
\begin{proof}
This follows from the $C$-$A$- (resp. $A$-$C$-) bimodule isomorphisms
\begin{align*}
C(G-, -) \otimes_B B(F-, -) &\simeq C(GF-, -) \\
B(-, F-) \otimes_B C(-, G-) &\simeq C(-, GF-)
\end{align*}
cf. \cite[Th.~4.47]{kelly}.
\end{proof}
\end{corollary}

\paragraph{Remarks}
\begin{enumerate}[(1)]
\item For right modules, the adjoints to $\mathrm{res}_F: \ModR{B} \to \ModR{A}$ (by which we technically mean $\mathrm{res}_{F^{\op}}$) are $F_L := B(-,F-) \otimes_{A^{\op}} \,?\, \simeq\, \,?\, \otimes_A B(-, F-)$ and $F_R := [B(F-,-), \,?\,]_{A^{\op}}$. The rest of the proposition remains the same except corresponding (c) and (d) statements, of (vi) and (vii), must be interchanged.
\item Note $F_L A(a, -) = B(F-,-) \otimes_A A(a, -) \simeq B(F(a), -)$ by the co-Yoneda isomorphisms, so $F_L$ preserves representable functors. $F_L$ is also right exact, and therefore restricts to a functor $\Rmod{A} \to \Rmod{B}$.
\end{enumerate}

\subsection{Tensor products and inverse limits}
If $M \in \Rmod{R}$, then $- \otimes_R M: \ModR{R} \to \Ab$ does not commute with all inverse limits. However, we have a partial result if we consider only countable inverse systems in $\modR{R}$.

\paragraph{}
An inverse system $((X_i)_{i \in \mathbb{N}}, (\gamma_{ji}: X_j \to X_i)_{j \ge i})$ of non-empty sets is a \textbf{Mittag-Leffler system} \cite[13.1.2]{grothendieck} if for all $i \ge 1$ there exists $j \ge i$ such that for all $k \ge j$ we have $\gamma_{ki}(X_k) = \gamma_{ji}(X_j)$. Equivalently, for all $i \ge 1$, the descending sequence
\begin{equation*}
\dots \subseteq \gamma_{i+2,i}(X_{i+2}) \subseteq \gamma_{i+1,i}(X_{i+1}) \subseteq X_i
\end{equation*}
eventually stabilizes.

\begin{lemma}
If $((X_i)_{i \in \mathbb{N}}, (\gamma_{ji}: X_j \to X_i)_{j \ge i})$ is a Mittag-Leffler system, then the inverse limit $\varprojlim_i X_i$ is non-empty.
\begin{proof}
For $i \in \mathbb{N}$ define $Y_i = \cap_{j \ge i} \gamma_{ji}(X_j) \subseteq X_i$, then for $j \ge i$ the function $\gamma_{ji}: X_j \to X_i$ restricts to a surjection $\gamma_{ji}: Y_j \to Y_i$. Moreover, the limits $\varprojlim X_i$ and $\varprojlim Y_i$ coincide, with the latter non-empty: choose $y_1 \in Y_1$ and inductively $y_{i+1} \in \gamma^{-1}_{i+1,i}(y_i)$ for $i \ge 1$, then $\bar{y} = (y_i)_{i \in \mathbb{N}} \in \varprojlim Y_i$.
\end{proof}
\end{lemma}

The following result is (essentially) \cite[13.2.2]{grothendieck}; we translate the proof here for convenience.
\begin{proposition}\label{prop:MittagLeffler}
Let $(f_i: X_i \to Y_i)_{i \in \mathbb{N}}$ be morphisms in $\modR{R}$ compatible with two inverse systems $(\alpha_{ji}: X_j \to X_i)_{j \ge i}$ and $(\beta_{ji}: Y_j \to Y_i)_{j \ge i}$ --- i.e. $f_j \alpha_{ji} = \beta_{ji} f_i$ for all $j \ge i$. Let $\bar{f} = (f_i)_{i \in \mathbb{N}}$ be the induced morphism $\varprojlim X_i \to \varprojlim Y_i$. Given $\bar{y} = (y_i)_{i \in \mathbb{N}} \in \varprojlim Y_i$, then $\bar{y} \in \mathrm{im}(\bar{f})$ if and only if $f^{-1}_i(y_i) \ne \emptyset$ for all $i \ge 1$.
\begin{proof}
Given $\bar{y} = (y_i)_{i \in \mathbb{N}} \in \varprojlim Y_i$ set $U_i = f^{-1}_i(y_i) \subseteq X_i$ and suppose $U_i \ne \emptyset$ for all $i \ge 1$. We have an induced inverse system $(\gamma_{ji}: U_j \to U_i)_{j \ge i}$ by restriction of $(\alpha_{ji}: X_j \to X_i)$. If there exists $\bar{x} \in \varprojlim U_i$, then $\bar{x} \in \varprojlim X_i$ and $f(\bar{x}) = \bar{y}$ as required. Thus, by the above lemma, it is enough to show that $(U_i)_{i \in \mathbb{N}}$ is a Mittag-Leffler system.

Let $L_i = \mathrm{ker}(f_i)$ for $i \in \mathbb{N}$, then $(L_i)_{i \in \mathbb{N}}$ is a Mittag-Leffler system (as each $L_i$ is a finite-dimensional $\field$-vector space). Given $i \in \mathbb{N}$ we can choose $j \ge i$ such that $\alpha_{ki}(L_k) = \alpha_{ji}(L_j)$ for all $k \ge j$. We claim $\alpha_{ki}(U_k) = \alpha_{ji}(U_j)$ for $k \ge j$. Let $u_i \in \alpha_{ji}(U_j)$ be given, so $u_i = \alpha_{ji}(u_j)$ for some $u_j \in U_j$. Choose any $u^\prime_k \in U_k$ and set $u^\prime_j = \alpha_{kj}(u^\prime_k)$ and $u^\prime_i = \alpha_{ji}(u^\prime_j)$. Then $f_j(u_j - u^\prime_j) = y_j - y_j = 0$ and $u_j-u^\prime_j \in L_j$. Hence $u_i-u^\prime_i \in \alpha_{ji}(L_j) = \alpha_{ki}(L_i)$ and there exists $x_k \in L_k$ such that $u_i-u^\prime_i = \alpha_{ki}(x_k)$. Now $f_k(u^\prime_k +x_k) = y_k$ (so $u^\prime_k+x_k \in U_k$) and $\alpha_{ki}(u^\prime_k + x_k) = u^\prime_i + (u_i-u^\prime_i) = u_i$ giving $u_i \in \alpha_{ki}(U_k)$. Therefore $\alpha_{ji}(U_j) \subseteq \alpha_{ki}(U_k)$ as claimed, proving one direction.

The other direction is trivial, for if $(y_i)_{i \in \mathbb{N}} = \bar{f}((x_i)_{i \in \mathbb{N}}) = (f_i(x_i))_{i \in \mathbb{N}}$, then $x_i \in f^{-1}(y_i)$ and so $f^{-1}(y_i) \ne \emptyset$ for all $i \ge 1$.
\end{proof}
\end{proposition}

\begin{proposition}\label{prop:TensInvs}
If $((M_i)_{i \in \mathbb{N}}, (\gamma_{ji}: M_j \to M_i)_{j \ge i})$ is an inverse system in $\modR{R}$, then $(\varprojlim M_i) \otimes_R X = \varprojlim(M_i \otimes_R X)$ for any $X \in \Rmod{R}$.
\begin{proof}
Take a presentation of $X$ as a finitely-presented module:
\begin{equation}\label{eq:Xfp}
\xymatrix{A^k \ar[r]^f & A^l \ar[r]^g & X \ar[r] & 0}
\end{equation}
As $- \otimes_R X$ is right exact, it induces the following commutative diagram, where all but the left-most column are necessarily exact.
\begin{equation*}
\xymatrix{
\varprojlim M^k_i \ar[d]^{\bar{f}=(f_i)} & \dots \ar[r] & M^k_2 \ar[d]^{f_2} \ar[r] & M^k_1 \ar[d]^{f_1} \\
\varprojlim M^l_i \ar[d]^{\bar{g}=(g_i)} & \dots \ar[r] & M^l_2 \ar[d]^{g_2} \ar[r] & M^l_1 \ar[d]^{g_1} \\
\varprojlim (M_i \otimes_R X) \ar[d] & \dots \ar[r] & M_2 \ar[d] \otimes_R X \ar[r] & M_1 \ar[d] \otimes_R X \\
0 & & 0 & 0
}
\end{equation*}
We show the left-most column is also exact. Surjectivity of $\bar{g}$ is by Proposition~\ref{prop:MittagLeffler} applied to the bottom two rows, since each $g_i$ is surjective. It is clear that $\bar{g}\bar{f} = 0$ since $g_i f_i = 0$ for all $i \ge 1$. Finally, if $\bar{m} = (m_i)_{i \in \mathbb{N}} \in \mathrm{ker}(\bar{g})$, then $m_i \in \mathrm{ker}(g_i) = \mathrm{im}(f_i)$ and $f^{-1}_i(m_i) \ne \emptyset$ for all $i \ge 1$. Thus, applying Proposition~\ref{prop:MittagLeffler} to the top two rows, we have $\bar{x} \in \mathrm{im}(\bar{f})$ as required.

Let $M = \varprojlim M_i$ in $\ModR{R}$, then applying $M\otimes_R -$ to \eqref{eq:Xfp} gives
\begin{equation*}
(\varprojlim M_i) \otimes_R X = \mathrm{cok}(1_M \otimes_R f) \simeq \mathrm{cok}(\bar{f}) = \varprojlim (M_i \otimes_R X)
\end{equation*}
as claimed.
\end{proof}
\end{proposition}
 
\newpage
\singlespacing

{\small
\textsc{University of Manchester, Manchester, UK}

\emph{Email address}: \texttt{michael.bushell@postgrad.manchester.ac.uk}
}

\end{document}